\documentclass[11pt]{article}
\usepackage{amssymb}
\usepackage{times}

\title{Complexit\'e des bor\'eliens \`a coupes d\'enombrables.\indent}
\author{Dominique LECOMTE}
\date{\it ~Fund. Math.~\rm 165 (2000), 139-174}

\newcommand{\Ana}{{\it\Sigma}^{1}_{1}}
\newcommand{\Ca}{{\it\Pi}^{1}_{1}}
\newcommand{\Borel}{{\it\Delta}^{1}_{1}}

\newcommand{\boraone}{{\bf\Sigma}^{0}_{1}}
\newcommand{\boratwo}{{\bf\Sigma}^{0}_{2}}

\newcommand{\boraxi}{{\bf\Sigma}^{0}_{\xi}}

\newcommand{\bortwo}{{\bf\Delta}^{0}_{2}}
\newcommand{\borone}{{\bf\Delta}^{0}_{1}}
\newcommand{\borthree}{{\bf\Delta}^{0}_{3}}

\newcommand{\bormone}{{\bf\Pi}^{0}_{1}}
\newcommand{\bormtwo}{{\bf\Pi}^{0}_{2}}

\newtheorem{thm} {Th\'eor\`eme}

\newtheorem{cor} [thm] {Corollaire}
\newtheorem{lem} [thm] {Lemme}

\begin{document}

\maketitle

\noindent {\footnotesize {\bf R\'esum\'e.} Nous donnons, pour chaque niveau de complexit\'e $\Gamma$, une caract\'erisation du type ``test d'Hurewicz" des bor\'eliens d'un produit de deux espaces polonais ayant toutes leurs coupes d\'enombrables ne pouvant pas \^etre rendus $\Gamma$ par changement des deux topologies polonaises. }\bigskip\smallskip

\noindent\bf {\Large 1 Introduction.}\rm\bigskip

 Cet article fait suite \`a une \'etude entam\'ee dans [Le1], [Le2] et 
[Le3]. Il a pour objet de r\'epondre \`a une conjecture faite dans ce dernier, et 
peut pour l'essentiel \^etre lu ind\'ependamment de ces articles. Cependant, la lecture 
pr\'ealable de ces articles peut \'eclairer plusieurs points techniques pr\'esents dans les 
arguments d\'evelopp\'es ici. Ces travaux se situent 
dans le cadre de la th\'eorie descriptive des ensembles. Je 
renvoie le lecteur \`a [Ku] pour les notions de base de th\'eorie descriptive 
classique et \`a [Mo] pour les notions de th\'eorie descriptive effective. 
 Pour d\'eterminer la complexit\'e exacte d'un bor\'elien, on est amen\'e \`a montrer qu'il n'est 
pas d'une classe de Baire donn\'ee. Le th\'eor\`eme d'Hurewicz, rappel\'e ci-dessous, donne 
un crit\`ere pour la classe des $G_\delta$ (cf [SR]) :\bigskip

\noindent\bf Th\'eor\`eme.\it\ Soient $X$ un espace polonais et $A$ un bor\'elien de $X$. Les 
conditions suivantes sont \'equivalentes :\smallskip

\noindent (a) Le bor\'elien $A$ n'est pas $\bormtwo$.\smallskip

\noindent (b) Il existe une injection continue $u:2^\omega\rightarrow X$ telle que 
$$u^{-1}(A)=\{\alpha\in2^\omega~/~\exists~n~~\forall~m\geq n~~\alpha (m)=0\}.$$\rm

 Cet exemple des suites nulles \`a partir d'un certain rang peut \^etre remplac\'e par 
n'importe quel ensemble infini d\'enombrable sans point isol\'e de $2^\omega$. Il est 
appel\'e ``test d'Hurewicz" pour la classe des $G_\delta$. Ce th\'eor\`eme a \'et\'e g\'en\'eralis\'e 
aux autres classes de Baire par A. Louveau et J. Saint Raymond (cf [Lo-SR]).\bigskip

 On s'int\'eresse ici \`a une hi\'erarchie analogue \`a celle de Baire, sauf 
qu'au lieu de partir des ouverts-ferm\'es d'un espace polonais de dimension 0, on part 
des produits de 
deux bor\'eliens, chacun d'entre eux \'etant inclus dans un espace polonais. L'analogie 
devient plus claire quand on sait qu'\'etant donn\'es un espace polonais $X$ et un 
bor\'elien $A$ de $X$, on peut trouver une topologie polonaise plus fine que la 
topologie initiale sur $X$ (topologie ayant donc les m\^emes bor\'eliens), de dimension 0, 
et qui rende $A$ ouvert-ferm\'e. Pour notre probl\`eme, le fait de travailler dans les 
espaces de dimension 0 n'est donc pas une restriction r\'eelle. La d\'efinition qui suit 
appara\^\i t alors naturelle.

\vfill\eject

\noindent\bf D\'efinition.\it\ Soient $X$ et $Y$ des espaces polonais, et $A$ un bor\'elien de 
$X \times Y$. Si $\Gamma$ est une classe de Baire, on dira que 
$A$ est $potentiellement\ dans\ \Gamma$ $($ce qu'on notera 
$A \in \mbox{pot}(\Gamma))$ s'il existe des topologies polonaises de dimension 0, 
$\sigma$ $($sur $X)$ et $\tau$ $($sur $Y)$, plus fines que les topologies 
initiales, telles que $A$, consid\'er\'e comme partie de $(X, \sigma)\times (Y, \tau)$, soit dans $\Gamma$.\rm\bigskip

 La motivation pour l'\'etude de ces classes de Wadge potentielles trouve son 
origine dans l'\'etude des relations d'\'equivalence bor\'eliennes, et plus pr\'ecis\'ement 
dans l'\'etude du pr\'e-ordre suivant sur la collection des relations d'\'equivalence 
bor\'eliennes d\'efinies sur un espace polonais : 
$$E\leq F~~\Leftrightarrow ~~\exists~f~\mbox{bor\'elienne}~~E=(f\times f)^{-1}(F).$$
A l'aide de la notion de classe de Wadge potentielle, A. Louveau montre dans 
[Lo3] que la collection des relations d'\'equivalence $\boraxi$ n'est pas 
co-finale, et il en d\'eduit qu'il n'existe pas de relation maximum pour $\leq$.\bigskip

 On cherche \`a \'etablir des r\'esultats analogues au th\'eor\`eme d'Hurewicz pour les classes 
de Baire potentielles. Le r\'esultat principal de cet article \'etablit l'analogue du 
th\'eor\`eme d'Hurewicz pour la classe des ensembles potentiellement $G_\delta$. Dans 
[Le3], il y a la\bigskip

\noindent\bf Conjecture.\it\ Il existe un bor\'elien $B$ de $2^\omega\times2^\omega$, tel que pour 
tous espaces polonais $X$ et $Y$, et pour tout bor\'elien $A$ de $X\times Y$ \`a coupes verticales 
d\'enombrables, on a l'\'equivalence entre les conditions suivantes :\smallskip

\noindent (a) Le bor\'elien $A$ n'est pas $\mbox{pot}(\bormtwo)$.\smallskip

\noindent (b) Il existe des fonctions continues $u : 2^\omega \rightarrow X$ et 
$v : 2^\omega \rightarrow Y$ telles que $\overline{B} \cap (u\times v)^{-1}(A) = B$.\rm\bigskip

 L'essentiel de cet article va consister \`a analyser les bor\'eliens non $\mbox{pot}(\bormtwo)$ 
ayant leurs coupes horizontales et verticales d\'enombrables pour arriver 
progressivement \`a montrer le\bigskip

\noindent\bf Th\'eor\`eme 7.\it\ Il existe un bor\'elien $B$ de $\omega^\omega\times\omega^\omega$, tel que pour 
tous espaces polonais $X$ et $Y$, et pour tout bor\'elien $A$ de $X\times Y$ dont les coupes 
horizontales et verticales sont d\'enombrables, on a l'\'equivalence entre les 
conditions suivantes :\smallskip

\noindent (a) Le bor\'elien $A$ n'est pas $\mbox{pot}(\bormtwo)$.\smallskip

\noindent (b) Il existe $u : \omega^\omega \rightarrow X$ et 
$v : \omega^\omega \rightarrow Y$, hom\'eomorphismes sur leurs images, tels que l'on ait 
$\overline{B}\cap (u\times v)^{-1}(A) = B$.\rm\bigskip

 Ce bor\'elien $B$ sera une ``version uniforme" du test d'Hurewicz, c'est-\`a-dire un 
ensemble dont toutes les coupes verticales sont infinies d\'enombrables sans point 
isol\'e. Plus pr\'ecis\'ement, $B$ sera r\'eunion d\'enombrable de graphes d'hom\'eomorphismes 
de domaine ouvert-ferm\'e. On verra aussi qu'essentiellement, dans tout bor\'elien ayant 
ses coupes horizontales et verticales d\'enombrables et n'\'etant pas $\mbox{pot}(\bormtwo)$, on 
peut trouver, \`a un changement de topologie pr\`es, une telle r\'eunion se r\'eduisant \`a $A$ 
au sens du th\'eor\`eme 7. On cherchera entre autres \`a r\'eduire de telles r\'eunions 
entre elles.\bigskip

 L'hypoth\`ese de d\'enombrabilit\'e des coupes dans le th\'eor\`eme 7 peut sembler moins 
naturelle que par exemple l'hypoth\`ese ``$A$ est $\mbox{pot}(\boratwo)$". Mais cette 
derni\`ere n'est pas suffisante. En effet, les bor\'eliens \`a coupes verticales (ou 
horizontales) d\'enombrables sont $\mbox{pot}(\boratwo)$ (voir [Lo2]). Nous montrons que le th\'eor\`eme 7 
devient faux si on suppose seulement $A$ \`a coupes verticales d\'enombrables, en 
utilisant l'injectivit\'e de $u$ et $v$.

\vfill\eject

 Pour terminer cette introduction, nous pla\c cons le th\'eor\`eme 7 dans un contexte 
plus g\'en\'eral. Il vient en effet compl\'eter l'\'etude des bor\'eliens \`a coupes 
verticales d\'enombrables commenc\'ee dans [Le3]. Je renvoie 
le lecteur \`a cet article pour les rappels concernant la hi\'erarchie de Wadge, qui 
affine celle de Baire. On peut montrer que les seules classes de Wadge non stables 
par passage au compl\'ementaire contenues dans $\bortwo=\boratwo\cap\bormtwo$ sont les 
diff\'erences transfinies d'ouverts. On peut d\'efinir sans probl\`eme la notion d'ensemble
potentiellement dans $\Gamma$, o\`u $\Gamma$ est une classe de Wadge, en utilisant la
m\^eme d\'efinition que pr\'ec\'edemment. L'analogue du th\'eor\`eme 7 pour les diff\'erences  transfinies
d'ouverts est montr\'e dans [Le3] (voir th\'eor\`emes 3.5 et 3.6).  A ceci pr\`es que
l'hypoth\`ese est moins forte (``$A$ est potentiellement
$\borthree$" au  lieu de ``$A$ a ses coupes horizontales et verticales d\'enombrables"), et 
que la conclusion est moins forte (on n'a pas l'injectivit\'e des fonctions de 
r\'eduction). Comme cons\'equence de ces r\'esultats, on obtient le r\'esultat de synth\`ese suivant :\bigskip 

\noindent\bf Corollaire 9.\it\ Soit $\Gamma$ une classe de Wadge non stable par passage au 
compl\'ementaire. Alors il 
existe un bor\'elien $B_\Gamma$ de $\omega^\omega\times\omega^\omega$ et un ferm\'e $F_\Gamma$ contenant 
$B_\Gamma$, tels que pour tous espaces polonais $X$ et $Y$, et pour tout bor\'elien $A$ 
de $X\times Y$ ayant ses coupes horizontales et verticales d\'enombrables, on a 
l'\'equivalence entre les conditions suivantes :\smallskip

\noindent (a) Le bor\'elien $A$ n'est pas $\mbox{pot}(\Gamma )$.\smallskip
 
\noindent (b) Il existe des fonctions continues $u : \omega^\omega\rightarrow X$ et 
$v : \omega^\omega\rightarrow Y$ telles que 
$F_\Gamma\cap (u\times v)^{-1}(A) = B_\Gamma$.\rm\bigskip

 Il est \`a noter que $B_\Gamma$ et $F_\Gamma$ vont \^etre donn\'es de mani\`ere explicite, et 
que $B_\Gamma$ a ses coupes horizontales et verticales d\'enombrables si $\Gamma\subseteq 
\bormtwo$, ce qui est le cas significatif. On a donc en particulier que $B_\Gamma\notin 
\mbox{pot}(\Gamma)$ si $\Gamma\subseteq 
\bormtwo$. D'autre part, si $\Gamma$ est auto-duale (c'est-\`a-dire si $\Gamma$ est 
stable par 
passage au compl\'ementaire), ne pas \^etre dans $\Gamma$ signifie ne pas \^etre dans l'une 
des deux classes non auto-duales succ\'edant \`a $\Gamma$ dans l'ordre de Wadge 
(l'inclusion des classes). L'\'etude des classes de Wadge auto-duales peut donc \^etre 
ramen\'ee \`a celle des classes de Wadge non auto-duales.\bigskip

\noindent\bf Question.\rm ~Un probl\`eme ouvert est de savoir si on peut 
supprimer l'hypoth\`ese ``$A$ a ses coupes horizontales et verticales d\'enombrables" dans le corollaire 9.\bigskip\smallskip

\noindent\bf {\Large 2 Analyse des bor\'eliens \`a coupes d\'enombrables n'\'etant pas $\mbox{pot}(\bormtwo)$.}\rm\bigskip

 La d\'efinition qui suit donne un sens pr\'ecis \`a l'expression ``version uniforme du test 
d'Hurewicz" \'evoqu\'ee dans l'introduction.\bigskip

\noindent\bf D\'efinition.\it\ On dira que $(Z,T,(g_{m,p})_{(m,p)\in\omega^2})$ est une 
$situation\ g\acute en\acute erale$ si\smallskip

\noindent (a) $Z$ et $T$ sont des espaces polonais parfaits de dimension $0$ non vides.\smallskip

\noindent (b) $g_{m,p}$ est un hom\'eomorphisme de domaine $D_{g_{m,p}}$ (respectivement d'image)
 ouvert-ferm\'e de $Z$ (respectivement de $T$).\smallskip
 
\noindent (c) Pour $m\in\omega$, $(D_{g_{m,p}})_{p\in\omega}$ est une suite de domaines deux \`a 
deux disjoints dont la r\'eunion est dense dans 
$Z$. On note $g_m$ la fonction obtenue par recollement des $g_{m,p}$, pour $p$ entier.\smallskip

\noindent (d) Il existe un $G_\delta$ dense $G(g)$ de $\bigcap_{m\in\omega} D_{g_m}$ tel que 
${g[x] := \{g_m(x)~/~m\in\omega\}}$ soit sans point isol\'e, pour tout  
$x$ de $G(g)$.\rm

\vfill\eject

 L'id\'ee va \^etre de chercher le bor\'elien $B$ du th\'eor\`eme 7 sous la forme 
$\bigcup_{m\in\omega} \mbox{Gr}(g_m\lceil G(g))$, o\`u $(Z,T,(g_{m,p})_{(m,p)\in\omega^2})$ est 
une situation g\'en\'erale. Et aussi de montrer que dans chaque bor\'elien $A$, dont toutes 
les coupes sont d\'enombrables et n'\'etant pas $\mbox{pot}(\bormtwo)$, on peut trouver, \`a un 
changement de topologie pr\`es, une telle r\'eunion se r\'eduisant \`a $A$ au sens du 
th\'eor\`eme 7. On va donc \^etre amen\'es \`a r\'eduire une situation g\'en\'erale \`a une autre. C'est 
l'objet du th\'eor\`eme 1 qui suit. Il se trouve que pour assurer l'existence d'une telle 
r\'eduction, il nous faut des conditions suppl\'ementaires, aussi bien au d\'epart qu'\`a 
l'arriv\'ee. D'o\`u les deux d\'efinitions qui suivent. Toutes les conditions suppl\'ementaires 
de ces d\'efinitions seront utilis\'ees dans la preuve du th\'eor\`eme 1 qui suit, \`a l'exception 
de la condition (b) d'une situation d'arriv\'ee, dont l'int\'er\^et appara\^\i tra plus tard.\bigskip

\noindent\bf D\'efinition.\it\ On dira que $(Z,T,(g_{m,p})_{(m,p)\in\omega^2})$ est une 
$situation\ d'arriv\acute ee$ si\smallskip

\noindent (a) $(Z,T,(g_{m,p})_{(m,p)\in\omega^2})$ est une situation g\'en\'erale.\smallskip

\noindent (b) Le diam\`etre du domaine et de l'image de $g_{m,p}$ vaut au  plus $2^{-\Delta (m,p)}$, o\`u 
$\Delta :\omega^2\rightarrow \omega$ est injective.\smallskip

\leftline{$\begin{array}{ll} \!\!\!\mbox{(c)}\!\!\!\!\!
& \mbox{Pour~tout~entier~non~nul}~k~\mbox{et~pour~toute~suite}~u~\mbox{dans}~
(\omega^2)^{2k},~\mbox{on~a~l'implication}\cr\cr
& \pmatrix{ 
& \exists~U\in\borone,~\emptyset\!\not=\! U\subseteq Z\cr 
& \mbox{et}\cr 
& \!\!\!\!\!\!\forall~x\in U~~g_{u(0)}^{-1}g_{u(1)}^{}...
g_{u(2k-2)}^{-1}g_{u(2k-1)}^{}(x)\! =\!x}\!\Rightarrow\exists~i\! <\! 2k-1~~u(i)\! =\! u(i+1).
\end{array}$}\rm\bigskip

\noindent\bf D\'efinition.\it\ On dira que $(F,(f_{n,p})_{(n,p)\in\omega^2})$ est une  
$situation\ de\ d\acute epart$ si\smallskip

\noindent (a) $(F,F,(f_{n,p})_{(n,p)\in\omega^2})$ est une situation g\'en\'erale. \smallskip

\noindent (b) $F = \Pi_{i\in\omega} A_i \subseteq \omega^\omega$, o\`u chaque $A_i\subseteq\omega$ est fini.\smallskip

\noindent (c) Pour $x$ dans $D_{f_n}$, on a $x <_{\mbox{lex}} f_n(x)$, et pour $n$, $p$ entiers, 
il existe un entier $q(n,p)$ tel que si $q>q(n,p)$ et $x$ est dans $D_{f_{n,p}}$, 
$f_{n,p}(x)(q) = x(q)$.\smallskip

\noindent (d) Pour $(x,y)$ dans $(\bigcap_{n\in\omega} D_{f_n}\times F)\cap \overline{
\bigcup_{n\in\omega} \mbox{Gr}(f_n)}$, il existe $(y_k)$ tendant vers $y$ telle que 
$(x,y_k)$ appartienne \`a $\bigcup_{n\in\omega} \mbox{Gr}(f_n)$ pour tout entier $k$, et 
$x\not= y$.\smallskip

\noindent (e) On d\'efinit des relations sur $\bigcup_{p\in\omega} 
(\Pi_{i< p}~A_i)$, $n$ \'etant entier, par 
$$s~{\cal R}\! _n~t~\Leftrightarrow~\vert s\vert =\vert t\vert ~\mbox{et}~s\leq_{{\mbox{lex}}} t~\mbox{et}~
(N_s\times N_t)\cap \mbox{Gr}(f_n)\not=\emptyset\mbox{,}$$
$$s~{\cal R}~t~\Leftrightarrow~\exists~n\in\omega~~s~{\cal R}\! _n~t.$$
 Si $s~{\cal R}~t$, on pose $\psi(s,t) := \mbox{min}\{n\in\omega~/~s~{\cal R}\! _n~t\}$. 
On demande\smallskip

\noindent (i) \mbox{[}$s~{\cal R}~s~\Rightarrow~\psi (s,s)=0$\mbox{]}, et \mbox{[}$s^\frown j~{\cal R}~t^\frown j~
\Rightarrow~\psi (s^\frown j,t^\frown j) = \psi (s,t)$\mbox{]}.\smallskip

 On d\'efinit ensuite la relation sym\'etrique engendr\'ee par $\cal R$ :
$$s~{\cal T}~t~\Leftrightarrow ~s~{\cal R}~t~\mbox{ou}~t~{\cal R}~s.$$
On dira que $u\in [\bigcup_{p\in\omega} (\Pi_{i< p}~A_i)]^{<\omega}\setminus\{
\emptyset\}$ est une ${\cal T}\mbox{-}cha{\hat\imath}ne$ si pour tout $j<\vert u\vert -1$, $u(j)~{\cal T}~u(j+1)$. 
On d\'efinit la relation d'\'equivalence engendr\'ee par $\cal R$ :
$$s~{\cal E}~t \Leftrightarrow \exists~u~~{\cal T}\mbox{-cha\^\i ne}~u(0)=
s~\mbox{et}~u(\vert u\vert -1)=t.$$
On demande\smallskip

\noindent (ii) Si $s~{\cal E}~t$, il existe une unique $\cal T$-cha\^\i ne $u$, sans 
r\'ep\'etition de termes, telle que $u(0)=s$ et aussi $u(\vert u\vert -1)=t$.\rm

\vfill\eject

\begin{thm} Soient $(F,(f_{n,p})_{(n,p)\in\omega^2})$ une 
situation de d\'epart et $(Z,T,(g_{m,p})_{(m,p)\in\omega^2})$ une situation 
d'arriv\'ee. Alors il existe des injections continues $u\! :\! F \rightarrow G(g)$ 
et $v\! :\! F\rightarrow T$ telles que\smallskip

\noindent (a) Pour $(x,y)$ dans $\bigcup_{n\in\omega} \mbox{Gr}(f_n)$, $v(y)$ est dans 
$g[u(x)]$.\smallskip

\noindent (b) Pour $(x,y)$ dans $(\bigcap_{n\in\omega} D_{f_n}\times F)\cap\overline{\bigcup_
{n\in\omega} \mbox{Gr}(f_n)}\setminus (\bigcup_{n\in\omega} \mbox{Gr}(f_n))$, $v(y)$ est dans 
$\overline{g[u(x)]}\setminus g[u(x)]$.\end{thm}

\noindent\bf D\'emonstration.\rm\ Le sch\'ema de la preuve est 
semblable \`a celui de la d\'emonstration du th\'eor\`eme 2.12 de [Le3] ; la condition (iv) ci-dessous apporte des complications. On va construire\bigskip

\noindent - Une suite $(U_s)_{s\in\bigcup_{p\in\omega}\Pi_{i<p}~A_i}$ d'ouverts-ferm\'es non vides 
de $Z$.\bigskip

\noindent - Une suite $(V_s)_{s\in\bigcup_{p\in\omega}\Pi_{i<p}~A_i}$ d'ouverts-ferm\'es non vides de $T$.\bigskip
 
\noindent - Une injection $\Phi : \bigcup_{p\in\omega}(\Pi_{i<p}~A_i)^2\rightarrow\omega^2$.\bigskip

 On pose, si $s~{\cal R}~t$, 
$$n(s,t) := \mbox{min}~\{n\in\omega\cup \{-1\}~/~\forall~n<i<\vert s\vert ~~s(i)=t(i)\}\mbox{,}$$
$$w(s,t) := \Phi (s\lceil (n(s,t)+1),t\lceil (n(s,t)+1)).$$ 
La premi\`ere coordonn\'ee de $w(s,t)$ sera not\'ee $w_0(s,t)$. 
Soit $(O_q)$ une suite d'ouverts de $Z$ telle que $G(g) := \bigcap_{q\in\omega} O_q$. 
On demande \`a ces objets de v\'erifier les conditions suivantes :
$$\begin{array}{ll} 
& (i)~~~U_{s^\frown i}\times V_{s^\frown i}\subseteq 
(U_s\cap\bigcap_{q\leq\vert s\vert } O_q)\times V_s \cr 
& (ii)~~\delta (U_s),~\delta (V_s)~\leq 2^{-\vert s\vert } \cr 
& (iii)~s~{\cal R}~t~\Rightarrow ~V_t = g_{w(s,t)} [U_s]~\mbox{et}~w_0(s,t)\geq\psi (s,t)\cr 
& (iv)~~(\vert s\vert =\vert t\vert ,~s <_{\mbox{lex}} t~\mbox{et}~r<\vert s\vert ~\mbox{v\'erifie}~\forall~q\leq r~~
s~\not \!{{\cal R}\! _q}~t)\Rightarrow (U_s\times V_t)\cap [\bigcup_{q\leq r} \mbox{Gr}(g_q)] = \emptyset\cr 
& (v)~~~U_{s^\frown i}\cap U_{s^\frown j} = V_{s^\frown i}\cap V_{s^\frown j} = \emptyset~\mbox{si}~i\not= j
\end{array}$$
$\bullet$ Admettons ceci r\'ealis\'e. On d\'efinit les injections continues $u$ et $v$ par 
les formules 
$$\{ u(\alpha )\} := \bigcap_{n\in\omega} U_{\alpha\lceil n}\mbox{,}~~~
\{ v(\alpha )\} := \bigcap_{n\in\omega} V_{\alpha\lceil n}.$$
 Si $y=f_n(x)$, pour tout entier $m$ on a $x\lceil m~{\cal R}~y\lceil m$, et il 
existe un entier $m_0$ tel que si $m\geq m_0$, $n(x\lceil m,y\lceil m) = 
n(x\lceil m_0,y\lceil m_0)$, puisque $x(q)$ et $y(q)$ co\"\i ncident \`a partir d'un 
certain rang (on utilise la condition (c) d'une situation de d\'epart). Posons 
$w_0 := w(x\lceil m_0,y\lceil m_0)$. Par (iii) on a $V_{y\lceil m} 
= g_{w_0}[U_{x\lceil m}]$ si $m\geq m_0$. D'o\`u 
$$g_{w_0}(u(x))\in g_{w_0}[\bigcap_{m\geq m_0} U_{x\lceil m}]\subseteq \bigcap_{m\geq m_0}
g_{w_0}[U_{x\lceil m}] = \bigcap_{m\geq m_0} V_{y\lceil m} = \{ v(y)\}$$
et $(u(x),v(y))\in \mbox{Gr}(g_{w_0})$.

\vfill\eject

 Si maintenant $(x,y)\in (\bigcap_{n\in\omega} D_{f_n}\times F)\cap\overline{\bigcup_
{n\in\omega} \mbox{Gr}(f_n)}\setminus (\bigcup_{n\in\omega} \mbox{Gr}(f_n))$, soit $(y_k)$ telle que 
$(x,y_k)\in \bigcup_{n\in\omega} \mbox{Gr}(f_n)$ et $y_k$ tende vers $y$ (on utilise la 
condition (d) d'une situation de d\'epart). Par ce qui pr\'ec\`ede, on a que ${(u(x), v(y_k))\in 
\bigcup_{m\in\omega} \mbox{Gr}(g_m)}$, d'o\`u ${v(y_k)\in g[u(x)]}$ et 
${v(y)\in\overline{g[u(x)]}}$, par continuit\'e de $v$. Pour tout entier $r$, il existe 
un entier ${p(r)>r}$ tel que pour $q\leq r$, 
$x\lceil p(r)~\!{\not\!{{\cal R}\! _q}}~y\lceil p(r)$ et $x\lceil p(r)<_{\mbox{lex}} 
y\lceil p(r)$ puisque $x\not= y$ ; par (iv) on a donc ${(U_{x\lceil p(r)}\times V_{y\lceil p(r)})\cap 
[\bigcup_{q\leq r} \mbox{Gr}(g_q)] = \emptyset}$ et 
$(u(x),v(y))\notin \bigcup_{q\leq r} \mbox{Gr}(g_q)$, donc $v(y)\notin g[u(x)]$.\bigskip

\noindent $\bullet$ Montrons donc que la construction est possible.\bigskip

 On pose $U_\emptyset := D_{g_{0,0}}$, 
$V_\emptyset := g_{0,0}[U_\emptyset]$, et $\Phi (\emptyset,\emptyset) := (0,0)$ (on peut 
supposer que $D_{g_{0,0}}\not=\emptyset$). Admettons avoir construit 
$U_s$, $V_s$ et $\Phi (s,t)$ pour $\vert s\vert ,~\vert t\vert \leq p$ v\'erifiant (i)-(v), et 
soient $s\in \Pi_{j<p}~A_j$ et $i\in A_p$. On note, si $\tilde s~{\cal E}~\tilde t$, 
${\cal T}(\tilde s,\tilde t)$ l'unique $\cal T$-cha\^\i ne sans r\'ep\'etition de termes 
$c$ telle que $c(0) = \tilde s$ et $c(\vert c\vert -1) = \tilde t$, qui existe par la 
condition (e).(ii) d'une situation de d\'epart.\bigskip

 La relation $\cal T$ d\'efinit sur ${\cal E}(s^\frown i)$ une structure d'arbre, de 
sommet $s^\frown i$. On va essentiellement construire les ouverts-ferm\'es cherch\'es par 
r\'ecurrence sur les niveaux de cet arbre. Si $k\in\omega$, on pose 
$H_k := \{z\in {\cal E}(s^\frown i)~/~\vert {\cal T}(z,s^\frown i)\vert = k+1\}$. Alors $H_k$ et le nombre de $H_k$ non vides sont finis, puisque les classes d'\'equivalence de $\cal E$ sont finies. De 
plus, $H_k$ est non vide si $H_{k+1}$ l'est, donc on peut trouver $q$ tel que $H_0$, ..., $H_q$ 
soient non vides et $H_k$ soit vide si $k>q$. Posons 
$$H_k := \{z_{(k,1)},...,z_{(k,p_k)}\}\mbox{,}~\phi : \left\{\!\!
\begin{array}{ll}
\bigcup_{k\leq q} \{ k\}\times\{1,...,p_k\}\!\!\!\! 
& \rightarrow \omega\cr 
(k,r) 
& \mapsto (\Sigma_{i<k}~p_i)+r
\end{array}
\right.$$
 On a donc Im$(\phi ) = \{p_0=1,p_0 +1,...,p_0 + p_1,...,p_0 +...+p_{q-1} +1, ..., 
p_0+...+p_q\}$. On va construire par r\'ecurrence sur $n\in\{1,...,p_0+...+p_q\}$, et 
pour $k\in \{1,...,n\}$, des ouverts-ferm\'es non vides $U^n_k$ et 
$V^n_k$ tels que, si $w(k,l) := w(z_{\phi^{-1}(k)},z_{\phi^{-1}(l)})$, on ait 
$$\begin{array}{ll} 
& (1)~U^n_k\times V^n_k\subseteq (U_{z_{\phi^{-1}(k)}\lceil p}\cap \bigcap_{q\leq p} O_q)\times 
V_{z_{\phi^{-1}(k)}\lceil p} \cr 
& (2)~\delta (U^n_k),~\delta (V^n_k)~\leq 2^{-p-1} \cr 
& (3)~\mbox{Si}~k,l\in\{ 1,...,n\}\mbox{, alors}~V^n_k = g_{0,0}[U^n_k] ~\mbox{et~si~de~plus}~
z_{\phi^{-1}(k)}~{\cal R}~z_{\phi^{-1}(l)}\mbox{,}\cr 
& \mbox{alors}~V^n_l = g_{w(k,l)}[U^n_k]~\mbox{et}~w_0(k,l)\geq \psi(z_{\phi^{-1}(k)},z_{\phi^{-1}(l)}) \cr 
& (4)~\mbox{Si}~k,l\in\{ 1,...,n\}\mbox{,}~z_{\phi^{-1}(k)} <_{\mbox{lex}} z_{\phi^{-1}(l)}~\mbox{et}~r \leq p~
\mbox{v\'erifie}~~\forall~q\leq r~~z_{\phi^{-1}(k)}~\! {\not \!{{\cal R}\! _q}}~z_{\phi^{-1}(l)}\mbox{,}\cr 
& \mbox{alors}~(U^n_k\times V^n_l)\cap[\bigcup_{q\leq r} \mbox{Gr}(g_q)] = \emptyset \cr 
& (5)~\mbox{Si}~k,l\in\{ 1,...,n\}~\mbox{et}~k\not= l\mbox{, alors}~U^n_k\cap U^n_l=V^n_k\cap V^n_l=\emptyset\cr 
& (6)~\mbox{Si}~k\in\{ 1,...,n\}\mbox{, alors}~U^{n+1}_k\times V^{n+1}_k \subseteq U^n_k\times V^n_k
\end{array}$$
Admettons cette construction effectu\'ee. On est tent\'e de poser 
$U_{z_{\phi^{-1}(k)}} := U^{p_0+...+p_q}_k$,  
$$V_{z_{\phi^{-1}(k)}} := V^{p_0+...+p_q}_k\mbox{,}$$ 
et si $\Phi (z,z')$ n'est pas encore d\'efini, on le d\'efinirait en assurant l'injectivit\'e de $\Phi$ et la seconde 
partie de la condition (iii) ; on obtiendrait alors (i), (ii), (iii), et aussi (iv), (v) dans ${\cal E}(s^\frown i)$. Mais il se pourrait qu'il y ait plusieurs classes d'\'equivalence dans $\Pi_{i\leq p}~A_i$. On remarque alors que les 
conditions (i), (ii), (iv) et (v) sont h\'er\'editaires. La construction montrera 
qu'on peut proc\'eder comme suit pour obtenir les conditions (iv) et (v) dans 
$\Pi_{i\leq p}~A_i$.

\vfill\eject

 Pour chaque couple $(s,t)$ de suites non \'equivalentes, on diminue 
les ouverts-ferm\'es concern\'es (\`a savoir $U_s$ et $V_t$ pour la condition (iv), 
($U_s$ et $U_t$) ou ($V_s$ et $V_t$) pour la condition (v)), de fa\c con \`a assurer la 
condition ((iv) ou (v)) pour ce couple (c'est possible par raret\'e des graphes). 
Puis on assure la condition (iii), ce qui diminue les ouverts-ferm\'es et permet de 
conserver les conditions (i) \`a (v) dans chaque classe. Comme les ensembles $A_i$ 
sont finis, on arrive ainsi \`a satisfaire les conditions (i) \`a (v) en un nombre fini 
d'\'etapes.\bigskip

\noindent $\bullet$ Montrons donc que cette nouvelle construction est possible. Si $n=1$, 
$\phi^{-1}(n)$ vaut $(0,1)$ et $z_{\phi^{-1}(n)} = s^\frown i$ ; on choisit pour 
$U^1_1$ un ouvert-ferm\'e de $(\bigcap_{q\leq p} O_q) \cap g_{0,0}^{-1}(V_s)$ de 
diam\`etre au plus $2^{-p-1}$ tel que $\delta (g_{0,0}[U^1_1]) \leq 2^{-p-1}$, et on 
pose $V^1_1 := g_{0,0}[U^1_1]$.\bigskip

 Admettons avoir construit les suites finies 
$U^1_1$, $V^1_1$, ..., $U^{n-1}_1$, 
$V^{n-1}_1$, ..., $U^{n-1}_{n-1}$, 
$V^{n-1}_{n-1}$, v\'erifiant (1)-(6), ce qui est fait pour $n=2$. On note 
${\cal T}(r,q) := {\cal T}(z_{\phi^{-1}(r)},z_{\phi^{-1}(q)})$. La suite 
$z_{\phi^{-1}(n)}$ est dans $H_{(\phi^{-1}(n))_0}$, donc on a 
$\vert {\cal T}(1,n)\vert -1 = (\phi^{-1}(n))_0$. Comme $p_0=1$, 
$(\phi^{-1}(n))_0\geq 1$, donc 
$\vert {\cal T}(1,n)\vert\geq 2$ et 
${\cal T}(1,n)(\vert {\cal T}(1,n)\vert -2)
\in H_{(\phi^{-1}(n))_0-1}$ ; par le choix de $\phi$, 
on peut trouver $m<n$ tel que ${\cal T}(1,n)(\vert {\cal T}
(1,n)\vert -2) = z_{\phi^{-1}(m)}$. D'o\`u 
$z_{\phi^{-1}(n)}~{\cal T}~z_{\phi^{-1}(m)}$. Notons 
$$o := \left\{\!\!\!\!\!\!
\begin{array}{ll} 
& n(z_{\phi^{-1}(n)},z_{\phi^{-1}(m)})~\mbox{si}~z_{\phi^{-1}(n)}~{\cal R}~z_{\phi^{-1}(m)}\mbox{,}\cr\cr 
& n(z_{\phi^{-1}(m)},z_{\phi^{-1}(n)})~\mbox{si}~z_{\phi^{-1}(m)}~{\cal R}~z_{\phi^{-1}(n)}.
\end{array}
\right.$$
\bf Cas 1.\rm\ $o<p$.\bigskip
 
\noindent 1.1. $z_{\phi^{-1}(m)}~{\cal R}~z_{\phi^{-1}(n)}$.\bigskip

 La suite $w(m,n) = \Phi(z_{\phi^{-1}(m)}\lceil (o+1),z_{\phi^{-1}(n)}\lceil (o+1))$ a 
d\'ej\`a \'et\'e d\'efinie et on a 
$$V_{z_{\phi^{-1}(n)\lceil p}} = g_{w(m,n)} [U_{z_{\phi^{-1}(m)\lceil p}}].$$
On choisit un ouvert-ferm\'e $\tilde U^n_n$ de 
$(\bigcap_{q\leq p} O_q)\cap g_{0,0}^{-1}(g_{w(m,n)}[U^{n-1}_m])$ de diam\`etre au plus $2^{-p-1}$ tel que $\delta (g_{0,0}[\tilde U^n_n])\leq 2^{-p-1}$, on pose 
$\tilde V^n_n := g_{0,0}[\tilde U^n_n]$. De sorte que les conditions (1) \`a (5) pour $k = l = n$ sont 
r\'ealis\'ees.\bigskip

 On d\'efinit ensuite les $\tilde U^n_q$ et $\tilde V^n_q$ pour 
$1\!\leq\! q\! <\! n$, par r\'ecurrence sur $\vert {\cal T}(q,n)\vert$. 
Comme $\vert {\cal T}(q,n)\vert\!\geq\! 2$, $\tilde U_{{\cal T}(q,n)(1)}^n$ a \'et\'e d\'efini et 
il y a 2 cas. Soit $r$ entier compris entre $1$ et $n$ tel que 
${\cal T}(q,n)(1) = z_{\phi^{-1}(r)}$ (la d\'efinition de 
$\phi$ montre l'existence de $r$).\bigskip

\noindent 1.1.1. $z_{\phi^{-1}(r)}~{\cal R}~z_{\phi^{-1}(q)}$.\bigskip

 On pose 
$$\left\{\!\!
\begin{array}{ll} 
\tilde V^n_q\!\!\!\! 
& := g_{w(r,q)}[\tilde U^n_r]\mbox{,}\cr  
\tilde U^n_q\!\!\!\! 
& := g_{0,0}^{-1}(\tilde V^n_q).
\end{array}
\right.$$
1.1.2. $z_{\phi^{-1}(q)}~{\cal R}~z_{\phi^{-1}(r)}$.\bigskip

 On pose 
$$\left\{\!\!
\begin{array}{ll}
\tilde U^n_q\!\!\!\! 
& := g^{-1}_{w(q,r)} (\tilde V^n_r)\mbox{,}\cr  
\tilde V^n_q\!\!\!\! 
& := g_{0,0}[\tilde U^n_q].
\end{array}\right.$$

\vfill\eject

Montrons que ces d\'efinitions sont licites. On a que 
${\cal T}(q,n)(1) = z_{\phi^{-1}(r)}$, o\`u 
$1\leq r \leq n$. Si le cas $r = n$ se produit, comme $z_{\phi^{-1}(m)}$ et 
$z_{\phi^{-1}(q)}$ sont dans ${\cal E}(s^\frown i)$, l'unicit\'e de 
${\cal T}(1,n)$ montre que $q = m$. On en d\'eduit que si 
$r = n$, on est dans le cas 1.1.2 puisqu'on ne peut pas avoir 
$z_{\phi^{-1}(q)}~=~z_{\phi^{-1}(r)}$, ${\cal T}(q,n)$ \'etant sans r\'ep\'etition de termes 
(si $\tilde s~{\cal R}~\tilde t$ et $\tilde t~{\cal R}~\tilde s$, 
on a $\tilde s = \tilde t$).\bigskip

 Dans le cas 1.1.1, on a $r<n$ et 
$V^{n-1}_q = g_{w(r,q)}[U^{n-1}_r]$, donc 
$\tilde V^{n}_q$ est un ouvert-ferm\'e non vide de 
$V^{n-1}_q$, puisque 
$\tilde U^n_r\subseteq U^{n-1}_r$. Par suite, $\tilde U^n_q$ est un ouvert-ferm\'e non 
vide de $U^{n-1}_q$. De m\^eme, 
$\tilde U^n_q$ est un ouvert-ferm\'e non vide de 
$U^{n-1}_q$ dans le cas 1.1.2, $r<n$. Si $r=n$, $q=m$ et la m\^eme 
conclusion vaut, par le choix de $\tilde U^{n}_n$ ; par suite, 
$\tilde V^n_q$ est un ouvert-ferm\'e non vide de 
$V^{n-1}_q$. D'o\`u la condition (6). Les conditions (1) et (2) pour $k=q$ en d\'ecoulent. V\'erifions (3). Soient donc $k,l\leq n$ tels que $z_{\phi^{-1}(k)}~{\cal R}~z_{\phi^{-1}(l)}$. Si on a l'\'egalit\'e
${\vert {\cal T}(k,n)\vert = \vert {\cal T}(l,n)\vert = 1}$, $k=l=n$, et la condition (3) est 
r\'ealis\'ee (on utilise la condition (e).(i) d'une situation de d\'epart). Plus g\'en\'eralement, 
la condition (3) est r\'ealis\'ee si $k=l$. Si 
$\vert {\cal T}(k,n)\vert = 1$ et  $\vert {\cal T}(l,n)\vert  = 2$, la liaison entre $z_{\phi^{-1}(k)}$ et 
$z_{\phi^{-1}(l)}$ a d\'ej\`a \'et\'e prise en compte, par minimalit\'e des longueurs. De m\^eme 
si ${\vert {\cal T}(k,n)\vert\! =\! 2}$ et $\vert {\cal T}(l,n)\vert  = 1$. Si $\vert {\cal T}(k,n)\vert$, $\vert {\cal T}(l,n)\vert\geq 2$, on a 
${\cal T} (k,n)(1)\! =\! {\cal T} (l,n)(0)$ ou bien alors ${\cal T} (k,n)(0)\! =\! 
{\cal T} (l,n)(1)$, par unicit\'e. L\`a encore, la liaison a 
\'et\'e prise en compte. On a donc $\tilde V_l^n = g_{w(k,l)}[\tilde U_k^n]$. La 
condition (3) est donc r\'ealis\'ee, le seul nouveau cas \'etant celui o\`u $r=n$ et $q=m$, et par 
la condition (e).(i) on a 
$$\begin{array}{ll}
w_0(z_{\phi^{-1}(m)},z_{\phi^{-1}(n)})\!\!\!\! 
& = w_0(z_{\phi^{-1}(m)}\lceil p,z_{\phi^{-1}(n)}\lceil p) \cr 
& \geq\psi(z_{\phi^{-1}(m)}\lceil p,z_{\phi^{-1}(n)}\lceil p) \cr 
& \geq\psi (z_{\phi^{-1}(m)},z_{\phi^{-1}(n)}).
\end{array}$$ 
Remarquons que les conditions (1), (2), (4), (5), (6) sont h\'er\'editaires. Soient 
$k,l\in\{ 1,...,n\}$ tels que $z_{\phi^{-1}(k)} <_{\mbox{lex}} z_{\phi^{-1}(l)}$ 
et $r\leq p$ tel que $\forall~q\leq r$, $z_{\phi^{-1}(k)}~{\not \!{{\cal R}\! _q}}~
z_{\phi^{-1}(l)}$. On veut assurer que 
$$(U^n_k\times V^n_l)\cap [\bigcup_{q\leq r} 
\mbox{Gr}(g_q)] =\emptyset .$$ 
On a $\tilde V^n_l = \tilde g_{k,l} [\tilde U^n_k]$, o\`u $\tilde g_{k,l}$ est de la forme 
$hh_{\vert {\cal T}(k,l)\vert-2}...h_0$, avec 
$$h := \left\{\!\!\!\!\!\!\!
\begin{array}{ll} 
& \mbox{Id}_T~\mbox{si}~\mbox{Im}(h_{\vert {\cal T}(k,l)\vert-2})\subseteq T\mbox{,}\cr 
& g_{0,0}~\mbox{sinon,}
\end{array}\right.$$
$$h_0 := \left\{\!\!\!\!\!\!\!
\begin{array}{ll} 
& g_{w({\cal T}(k,l)(0),{\cal T}(k,l)(1))}~\mbox{si}~{\cal T}(k,l)(0)~{\cal R}~{\cal T}(k,l)(1)\mbox{,}\cr 
& g_{w({\cal T}(k,l)(1),{\cal T}(k,l)(0))}^{-1}g^{}_{0,0}~\mbox{sinon,}
\end{array}\right.$$
$$h_{p+1} := \left\{\!\!\!\!\!\!\!
\begin{array}{ll} 
& g_{w({\cal T}(k,l)(p+1),{\cal T}(k,l)(p+2))}~\mbox{si}~{\cal T}(k,l)(p+1)~{\cal R}~{\cal T}(k,l)(p+2)~\mbox{et}~\mbox{Im}(h_p)\subseteq Z\mbox{,}\cr 
& g_{w({\cal T}(k,l)(p+2),{\cal T}(k,l)(p+1))}^{-1}g^{}_{0,0}~\mbox{si}~{\cal T}(k,l)(p+2)~{\cal R}~
{\cal T}(k,l)(p+1)~\mbox{et}~\mbox{Im}(h_p)\subseteq Z\mbox{,}\cr 
& g^{}_{w({\cal T}(k,l)(p+1),{\cal T}(k,l)(p+2))}g_{0,0}^{-1}~\mbox{si}~{\cal T}(k,l)(p+1)~{\cal R}~
{\cal T}(k,l)(p+2)~\mbox{et}~\mbox{Im}(h_p)\subseteq T\mbox{,}\cr 
& g_{w({\cal T}(k,l)(p+2),{\cal T}(k,l)(p+1))}^{-1}~\mbox{si}~{\cal T}(k,l)(p+2)~{\cal R}~{\cal T}(k,l)(p+1)~
\mbox{et}~\mbox{Im}(h_p)\subseteq T.
\end{array}\right.$$
En renum\'erotant, on a donc que $\tilde g_{k,l} = g_{s_{2r}}^{}g_{s_{2r-1}}^{-1}...g_{s_{1}}^{-1}g_{s_{0}}^{}$. On va donc chercher \`a assurer que 
$$(\tilde U^n_k\times \tilde g_{k,l} [\tilde U^n_k])\cap [\bigcup_{q\leq r} \mbox{Gr}(g_q)] =\emptyset .$$

\vfill\eject

 Montrons par l'absurde que c'est possible. Comme les $g_{m,p}$ sont des 
hom\'eomorphismes, $\tilde g_{k,l}$ est d\'efinie et continue donc on peut trouver un ouvert-ferm\'e 
non vide $U\subseteq \tilde U^n_k$, $q\leq r$ et $j\in \omega$ tels que 
pour tout $x\in U$, $g_q (x) = g_{q,j}(x)$ et $\tilde g_{k,l}(x) = g_q(x)$. En d'autres 
termes, $g^{-1}_{q,j}\tilde g^{}_{k,l}(x)=x$. Comme 
$(Z,T,(g_{m,p})_{(m,p)\in\omega^2})$ est une situation d'arriv\'ee, 
on a ${ s_{2r}=(q,j)}$ ou alors il existe $i<2r$ tel que $s_i = s_{i+1}$. Montrons que 
la seconde \'eventualit\'e est exclue. Deux types de cas peuvent se produire :
$$\begin{array}{ll} s_i = w({\cal T}(k,l)(j),{\cal T}(k,l)(j+1))\!\!\!\!
& \mbox{et}\ \ \ \ s_{i+1} = w({\cal T}(k,l)(j+2),{\cal T}(k,l)(j+1)) \cr 
& \mbox{ou}\cr 
\ \ \ \ \ \ \ \ \ \ \ \ \ \ \ \ \ \ \ \ \ \ \ \ \ \ \ \ \ \ \ \ \ \ \ \ \ \ \ s_i = (0,0)
& \mbox{et}\ \ \ \ s_{i+1} = w({\cal T}(k,l)(j),{\cal T}(k,l)(j+1)).
\end{array}$$
La fonction $\Phi$ \'etant injective, ${\cal T}(k,l)(j)\lceil (n({\cal T}(k,l)(j),{\cal T}(k,l)(j\! +\! 1))\! +\! 1)$ vaut 
$${\cal T}(k,l)(j\! +\! 2)\lceil (n({\cal T}(k,l)(j\! +\! 2),{\cal T}(k,l)(j\! +\! 1))\! +\! 1)$$ 
et $n({\cal T}(k,l)(j),{\cal T}(k,l)(j+1))+1 = n({\cal T}(k,l)(j+2),{\cal T}(k,l)(j+1))+1$. Par d\'efinition de $n(s,t)$, on a ${{\cal T}(k,l)(j) = {\cal T}(k,l)(j+2)}$, ce qui contredit la 
d\'efinition de ${\cal T}(k,l)$.\bigskip

 Dans l'autre type de cas, l'injectivit\'e de $\Phi$ et la d\'efinition de $n(s,t)$ font que 
$${{\cal T}(k,l)(j) = {\cal T}(k,l)(j+1)}\mbox{,}$$ 
ce qui est \'egalement absurde.\bigskip

 On a donc 
${s_{2r}= (q,j)}$, $r=0$ et $z_{\phi^{-1}(k)}~{\cal R}~z_{\phi^{-1}(l)}$ ;   
${w_0({\cal T}(k,l)(0),{\cal T}(k,l)(1))\! =\! q\!\leq\! r}$, et  
${\psi(z_{\phi^{-1}(k)},z_{\phi^{-1}(l)})\leq w_0
(z_{\phi^{-1}(k)},z_{\phi^{-1}(l)})\leq r}$ ; d'o\`u   
$z_{\phi^{-1}(k)}~{\cal R}\! _{\psi(z_{\phi^{-1}(k)},z_{\phi^{-1}(l)})}~z_{\phi^{-1}(l)}$, avec 
$$\psi(z_{\phi^{-1}(k)},z_{\phi^{-1}(l)})\leq r\mbox{,}$$
ce qui est absurde.\bigskip

 En diminuant $\tilde U_k^n$ et en assurant la condition (3), on peut donc avoir la 
condition (4) pour le couple $(k,l)$. Comme ces couples sont en nombre fini, on 
obtient la condition (4) en un nombre fini d'\'etapes.\bigskip

 On doit maintenant assurer (5). Il existe une relation du type 
${\tilde U^n_l = \tilde g_{k,l} [\tilde U^n_k]}$. On montre alors comme avant qu'on peut diminuer 
$\tilde U^n_k$ de fa\c con \`a assurer la disjonction de $\tilde U^n_k$ et $\tilde U^n_l$. 
De m\^eme, il existe une relation du type 
$\tilde V^n_l = \tilde g_{k,l} [\tilde V^n_k]$. L\`a encore, on peut diminuer 
$\tilde V^n_k$ 
de fa\c con \`a assurer la disjonction de $\tilde V^n_k$ et $\tilde V^n_l$. En effet, on a 
${\tilde g'_{k,l}(x) = x}$ pour tout $x$ d'un ouvert-ferm\'e non vide de $Z$, o\`u 
$\tilde g'_{k,l} :=  g_{(0,0)}^{-1}\tilde g_{k,l}g_{(0,0)}^{}$. On raisonne alors 
comme dans le cas pr\'ec\'edent pour avoir une contradiction.\bigskip

\noindent 1.2. $z_{\phi^{-1}(n)}~{\cal R}~z_{\phi^{-1}(m)}$.\bigskip

 Ce cas est analogue au pr\'ec\'edent (on a 
$V_{z_{\phi^{-1}(m)\lceil p}} = g_{w(n,m)} [U_{z_{\phi^{-1}(n)\lceil p}}]$, on choisit 
$\tilde U_n^n$ dans $g_{w(n,m)}^{-1}(V_m^{n-1})\cap(\bigcap_{q\leq p} O_q)$, et seul le 
cas 1.2.1 est possible si $r=n$).

\vfill\eject

\noindent\bf Cas 2.\rm\ $o=p$.\bigskip

\noindent 2.1. $z_{\phi^{-1}(m)}~{\cal R}~z_{\phi^{-1}(n)}$.\bigskip

 Alors $n(z_{\phi^{-1}(m)}\lceil\!  p,z_{\phi^{-1}(n)}\lceil\! p)\! <\! o\! =\! p$. On a 
$V_{z_{\phi^{-1}(n)}\lceil p}\! =\! g_{w(z_{\phi^{-1}(m)}\lceil p,z_{\phi^{-1}(n)}
\lceil p)}[U_{z_{\phi^{-1}(m)}\lceil p}]$, et $(U^{n-1}_m\times 
g_{w(z_{\phi^{-1}(m)}\lceil p,z_{\phi^{-1}(n)}\lceil p)}[U^{n-1}_m])\cap 
[\bigcup_{m\in\omega} \mbox{Gr}(g_m)]\not=\emptyset$. Par la condition (d) d'une situation 
g\'en\'erale, on peut trouver $q\geq \psi(z_{\phi^{-1}(m)},z_{\phi^{-1}(n)})$ et 
$j\in\omega$ tels que 
$${(U^{n-1}_m\times g_{w(z_{\phi^{-1}(m)}\lceil p,z_{\phi^{-1}(n)}\lceil p)}[U^{n-1}_m])
\cap \mbox{Gr}(g_{q,j})\not=\emptyset}.$$ 
On pose ${\Phi(z_{\phi^{-1}(m)},z_{\phi^{-1}(n)}) := (q,j)}$, en ayant pris soin de choisir $q$ 
suffisamment grand pour assurer l'injectivit\'e de $\Phi$.\bigskip

 L'ensemble $(\bigcap_{q\leq p} O_q)\cap g_{(0,0)}^{-1}(g_{w(m,n)} [U^{n-1}_m\cap 
g_{w(m,n)}^{-1}(g_{w(z_{\phi^{-1}(m)}\lceil p,z_{\phi^{-1}(n)}\lceil p)}
[U^{n-1}_m])])$ est un 
ouvert non vide ; on choisit 
$\tilde U^n_n$ dans cet ouvert et on raisonne comme en 1.1.1.\bigskip

\noindent 2.2. $z_{\phi^{-1}(n)}~{\cal R}~z_{\phi^{-1}(m)}$.\bigskip

 On choisit $\tilde U^n_n$ dans 
${g_{q,j}^{-1}(V^{n-1}_m)\cap g_{w(z_{\phi^{-1}(n)}\lceil p,
z_{\phi^{-1}(m)}\lceil p)}^{-1}(V^{n-1}_m)}\cap (\bigcap_{q\leq p} O_q)$ et on 
raisonne comme en 2.1, en posant 
$\Phi (z_{\phi^{-1}(n)},z_{\phi^{-1}(m)}) := (q,j)$.$\hfill\square$\bigskip

 Le lemme 6 de la section 3 assurera l'existence d'une situation de d\'epart. Le reste de 
cette section est consacr\'e \`a la recherche, dans chaque bor\'elien $A$ dont toutes les 
coupes sont d\'enombrables et n'\'etant pas $\mbox{pot}(\bormtwo)$, d'une r\'eunion 
$\bigcup_{m\in\omega} \mbox{Gr}(g_m\lceil G(g))$ se r\'eduisant \`a $A$, o\`u 
$(Z,T,(g_{m,p})_{(m,p)\in\omega^2})$ est une situation d'arriv\'ee. Pour ce faire, il se 
trouve qu'un r\'esultat interm\'ediaire (le th\'eor\`eme 3) va \^etre utilis\'e deux fois. Pour 
l'\'enoncer, il nous faut une d\'efinition suppl\'ementaire.\bigskip

\noindent\bf D\'efinition.\it\ On dira que $(Z,T,(h_{n,p})_{(n,p)\in\omega^2},M)$ est un 
$syst\grave eme\ r\acute educteur$ si\smallskip

\noindent (a) $Z$ et $T$ sont des espaces polonais parfaits de dimension $0$ non vides.\smallskip

\noindent (b) $h_{n,p}$ est un hom\'eomorphisme de domaine (respectivement d'image) ouvert-ferm\'e de 
$Z$ (respectivement de $T$), de diam\`etre au 
plus $2^{-\Delta (n,p)}$, o\`u $\Delta :\omega^2\rightarrow \omega$ est injective.\smallskip

\noindent (c) $M$ est $G_\delta$ de $Z\times T$ de projection dense dans $Z$, et 
$M\cap (\bigcup_{(n,p)\in\omega^2} \mbox{Gr}(h_{n,p}))=\emptyset$.\smallskip

\noindent (d) Pour tout ouvert $O$ de $Z\times T$ tel que $\Pi_Z[M\cap O]$ soit dense 
dans $Z$, $\Pi_Z[M\cap O]$ est comaigre dans $Z$.\smallskip

\noindent (e) Si $U$ et $V$ sont ouverts-ferm\'es et $M\cap (U\times V)\not=\emptyset$, 
$\mbox{Gr}(h_{n,p})\cap\overline{M}\cap (U\! \times\!  V)$ est le graphe de la restriction 
de $h_{n,p}$ \`a un ouvert $U_{n,p,U,V}$ de $U$ et ${\{ (n,p)\! \in\! 
\omega ^2 / \mbox{Gr}(h_{n,p})\! \cap\! \overline{M}\! \cap \! (U\! \times \! V)\! \not=\! 
\emptyset\}}$ est infini.\rm\bigskip

 Avec le lemme 2 et le th\'eor\`eme 3 qui suivent, nous reprenons pour l'essentiel la 
preuve du th\'eor\`eme 2.11 de [Le2] ; seul le vocabulaire change. Le lemme 2 
donne un proc\'ed\'e pour obtenir des syst\`emes r\'educteurs. Le lemme 5 nous en fournira un 
autre.

\vfill\eject

\begin{lem} Soient $X$, $Y$ des espaces polonais, $A$ un bor\'elien de $X\times Y$ 
dont les coupes horizontales et verticales sont d\'enombrables n'\'etant pas 
$\mbox{pot}(\bormtwo)$. Alors il existe un syst\`eme r\'educteur 
$(Z,T,(h_{n,p})_{(n,p)\in\omega^2},M)$ et des injections 
continues $u:Z\rightarrow X$ et $v:T\rightarrow Y$ tels que\smallskip

\noindent (a) $\bigcup_{(n,p)\in\omega^2} \mbox{Gr}(h_{n,p})\subseteq (u\times v)^{-1}(A)$.\smallskip

\noindent (b) $M\subseteq (u\times v)^{-1}(\check A)$.\end{lem}

\noindent\bf D\'emonstration.\rm\ On peut supposer, pour simplifier l'\'ecriture, que 
$X$ et $Y$ sont r\'ecursivement pr\'esen-t\'es, et que $A$ est 
$\Borel$-r\'eunion de graphes $\Borel$ de fonctions injectives. En effet, $A$ est 
r\'eunion d\'enombra-ble de graphes de fonctions bor\'eliennes, $A$ \'etant \`a coupes verticales 
d\'enombrables (cf [Ke]). $A$ \'etant \`a coupes horizontales d\'enombrables, ces 
fonctions sont countable-to-one ; on applique alors le lemme 2.4.(a) de [Le2] pour 
voir que $A$ est r\'eunion d\'enombrable de graphes d'injections bor\'eliennes. D\'esi-gnons 
par $W^X$ un 
ensemble $\Ca \subseteq \omega$ de codes pour les $\Borel$ de $X$, et 
par $C^X \subseteq \omega\times X$ un ensemble $\Ca$ dont les sections 
aux points de $W^X$ d\'ecrivent les $\Borel$ de $X$, et tel que la relation 
$$n\in W^X~\mbox{et}~(n,x)\notin C^X$$ soit $\Ca$ (cf [Lo1]). Soit \'egalement 
$W\subseteq W^{X\times Y}$ un ensemble $\Ca$ de codes pour les $\Borel 
\cap \mbox{pot}(\boratwo)$ de $X\times Y$ (dont l'existence 
est d\'emontr\'ee dans [Lo1]). Posons 
 $$H := \cup \{ (E\times F)\setminus A~/~E,F\in \Borel~\mbox{et}~(E\times F)
\setminus A\in \mbox{pot}(\boratwo) \}.$$
On a 
$$\begin{array}{ll}
H(x,y)~\Leftrightarrow\!\! & \exists~n\in W~\exists~m~(m)_0 \in W^X ~\mbox{et}~(m)_1 \in W^Y~\mbox{et}~\forall~z~~\forall~t \cr 
& [(n\in W^{X\times Y}~\mbox{et}~(n,z,t)\notin C^{X\times Y})~\mbox{ou}\cr 
& \{((m)_0,z)\in C^X~\mbox{et}~((m)_1,t)\in C^Y~\mbox{et}~(z,t)\notin A\}]~\mbox{et}\cr 
& [((m)_0\in W^X~\mbox{et}~((m)_0,z)\notin C^X)~\mbox{ou}~((m)_1\in W^Y~\mbox{et}~((m)_1,t)\notin C^Y)~\mbox{ou}\cr 
& (z,t)\in A~\mbox{ou}~(n,z,t)\in C^{X\times Y}]~\mbox{et}\cr 
& ((m)_0,x)\in C^X~\mbox{et}~((m)_1,y)\in C^Y~\mbox{et}~(x,y)\notin A.
\end{array}$$
Donc $H$ est $\Ca$. Posons $N := \check A \cap \check H$ ; $N$ est $\Ana$ et 
$G_{\delta}$ de $X\times Y$ muni de la topologie 
${\it\Delta}_X \times {\it\Delta}_Y$ (${\it\Delta}_X$ est la topologie 
engendr\'ee par les $\Borel$ de $X$). En effet, comme $A$ est $\mbox{pot}(\boratwo )$ et 
$\Borel$, $A$ est $\boratwo$ pour la topologie ${\it\Delta}_X \times {\it\Delta}_Y$ (voir 
[Lo1]). De m\^eme, $H$ est r\'eunion d\'enombrable de $\boratwo$ pour ${\it\Delta}_X \times {\it\Delta}_Y$. 
Posons $D_X := \{ x\in X~/~x\notin \Borel \}$, 
$\Omega_X := \{ x\in X~/~\omega_1^x = \omega_1^{\mbox{CK}} \}$, et 
$Z_0 := \Omega_X \cap D_X$, $T := \Omega_Y \cap D_Y$. On munit $Z_0$ 
(respectivement $T$) de la restriction de la topologie ${\it\Sigma}_X$ (respectivement 
${\it\Sigma}_Y$) de Gandy-Harrington sur $X$ (respectivement $Y$), de sorte que $Z_0$ et 
$T$ sont polonais parfaits de dimension 0. En effet, les traces des 
$\Ana$ sur $\Omega_X$ sont ouverts-ferm\'es de $(\Omega_X,{\it\Sigma}_X \lceil \Omega_X)$ : si $B$ 
est $\Ana$ contenu dans $\Omega_X$ et $f$ est $\Borel$ telle que
$f(x) \in WO \Leftrightarrow x \notin B$, on a 
$$x \notin B \Leftrightarrow x \notin \Omega_X~\mbox{ou}~(x\in \Omega_X~\mbox{et}~\exists~\xi <
\omega_1^{\mbox{CK}}~~(f(x) \in WO~\mbox{et}~ \vert f(x)\vert \leq \xi)).$$ 
L'espace $(\Omega_X,{\it\Sigma}_X \lceil \Omega_X)$ est donc \`a base d\'enombrable
d'ouverts-ferm\'es, donc m\'etrisable s\'eparable ; on sait (cf [HKL] et [Mo]) 
que c'est un espace fortement $\alpha$-favorable, comme ouvert (puisque $\Ana$) d'un 
espace fortement $\alpha$-favorable. C'est donc un espace polonais (cf [Ke]), 
de dimension 0 par ce qui pr\'ec\`ede.

\vfill\eject

 On pose $M := 
(Z_0\times T)\cap N$, $Z:=\Pi_{Z_0}[M]$. L'ensemble $M$ est non vide. En effet, $N$ 
rencontre $D_X\times D_Y$, sinon $N$ serait $\mbox{pot}(\borone)$ et 
$A\cup H$ aussi, donc $A = (A\cup H)\cap \check H$ serait 
$\mbox{pot}(\bormtwo)$, ce qui est exclus. Donc $N\cap (D_X\times D_Y)$, 
qui est $\Ana$, rencontre $\Omega_{X\times Y}$, par le th\'eor\`eme de base de Gandy 
(cf [Lo3]). Comme $\Omega_{X\times Y}\subseteq \Omega_X 
\times \Omega_Y$, $M\not= \emptyset$. Donc $Z$ et $T$ sont non vides et on a la 
condition (a) d'un syst\`eme r\'educteur, puisque $Z$ est $\Ana$, donc ouvert-ferm\'e de $Z_0$.
 De plus, $M$ est $G_\delta$ de $Z\times T$ et $\Pi_Z[M]=Z$ est dense dans $Z$.\bigskip
 
\noindent $\bullet$ Si $O$ est ouvert de $Z\times T$, il est r\'eunion de rectangles $\Ana$, donc 
$\Pi_Z[M\cap O]$ est ouvert de $Z$ ; par suite, si $\Pi_Z[M\cap O]$ est dense dans $Z$, 
$\Pi_Z[M\cap O]$ est comaigre dans $Z$. D'o\`u la condition (d) d'un syst\`eme r\'educteur.\bigskip

\noindent $\bullet$ Posons $A = \bigcup_{n\in\omega,~\mbox{disj.}} \mbox{Gr}(g_n)$, o\`u $g_n$ est un 
hom\'eomorphisme de domaine et d'image ouverts-ferm\'es non vides (en fait $\Borel$ ; le
raffinement des topologies permet de rendre les bijections bor\'eliennes bicontinues). Soient
$U$ et $V$ tels que ${M\! \cap\!  (U\! \times\!  V)\! \not=\!  \emptyset}$, ouverts-ferm\'es. 
Pour voir que 
$${\{n\in\omega~/~\mbox{Gr}(g_n)\cap\overline{M}\cap (U\times V)\not=\emptyset\}}$$ 
est infini, on peut supposer que $U$ et $V$ sont $\Ana$. Admettons avoir trouv\'e $n_1$, ..., $n_r$ deux \`a 
deux distincts tels que ${\mbox{Gr}(g_{n_i})\cap\overline{M}\cap (U\times V)\not=\emptyset}$
 pour $0<i\leq r$. Posons ${A' := \bigcup_{n\notin\{ n_1,...,n_r\}} \mbox{Gr}(g_n)}$. $A'$ 
est $\Borel$, et ${A=A'\cup\bigcup_{0<i\leq r} \mbox{Gr}(g_{n_i})}$. Posons 
${\cal O} := N\cap (U\times V)$ (${\cal O} = M\cap (U\times V)$ car ${U\times V
\subseteq Z\times T\subseteq Z_0\times T}$). Supposons que $A'\cap 
\overline{{\cal O}}^{{\it\Sigma}_X\times {\it\Sigma}_Y} = \emptyset$. Alors, par double application du th\'eor\`eme de s\'eparation, 
$A'\cap \overline{{\cal O}}^{{\it\Delta}_X\times {\it\Delta}_Y} = \emptyset$, donc on a la triple 
inclusion 
$${(U\times V)\setminus A'\subseteq {\cal O} \cup H\cup \bigcup_{0<i\leq r} 
\mbox{Gr}(g_{n_i}) \subseteq \overline{{\cal O}}^{{\it\Delta}_X\times {\it\Delta}_Y}\cup H \cup 
\bigcup_{0<i\leq r} \mbox{Gr}(g_{n_i})\subseteq \check A'}.$$
Donc $(U\times V)\setminus A'$ et $A'$ sont deux $\Ana$ s\'eparables par un 
$\mbox{pot}(\boratwo)$ ; ils peuvent par cons\'equent \^etre s\'epar\'es par un ensemble 
$K\in \Borel \cap \mbox{pot}(\boratwo)$ (cf [Lo1]). On a $U\times V \subseteq 
K\cup A'$, donc on peut trouver $\cal U$ et $\cal V$, deux $\Borel$ tels que 
$U\times V\subseteq {\cal U}\times {\cal V} \subseteq K\cup A'$. D'o\`u 
$({\cal U}\times {\cal V})\setminus A' = K\cap ({\cal U}\times {\cal V})$ et 
$({\cal U}\times {\cal V})\setminus A = K\cap ({\cal U}\times {\cal V})\cap 
[\bigcap_{0<i\leq r} ({\cal U}\times {\cal V})\setminus \mbox{Gr}(g_{n_i})]$ est 
$\mbox{pot}(\boratwo)$. On a donc $({\cal U}\times {\cal V})\setminus A\subseteq H$, puis 
${\cal O} \subseteq H\setminus H = \emptyset$, ce qui est absurde. On a donc que 
$A'\cap \overline{{\cal O}}^{{\it\Sigma}_X\times {\it\Sigma}_Y} \not= \emptyset$. Or ${\cal O}
\subseteq D_X\times D_Y$, donc $\overline{{\cal O}}^{{\it\Sigma}_X\times {\it\Sigma}_Y} 
\subseteq D_X\times D_Y$. Donc $(D_X\times D_Y)\cap A'\cap \overline{{\cal O}}^{{\it\Sigma}_X
\times {\it\Sigma}_Y} \not= \emptyset$ et est $\Ana$, donc rencontre 
$$\Omega_{X\times Y}\subseteq \Omega_X \times \Omega_Y\mbox{,}$$ 
donc $(Z_0\cap T)\cap A'\cap 
\overline{{\cal O}}^{{\it\Sigma}_X\times {\it\Sigma}_Y} \not= \emptyset$. Comme $U$ est 
ouvert-ferm\'e de $Z_0$, $A'\cap \overline{{\cal O}}^{Z\times T}\not= \emptyset$. L'existence 
d'un ouvert $U_{n,U,V}$ de $U$ tel que $\mbox{Gr}(g_n\lceil U_{n,U,V}) = 
\mbox{Gr}(g_n)\cap\overline{M}\cap (U\times V)$ vient du fait que $\overline{M}\in\Ana$.\bigskip

\noindent $\bullet$ Soit $\Delta : \omega^2\rightarrow \omega$ injective. On 
construit des ouverts-ferm\'es $D_{n,p}$ de $Z$ tels que\bigskip

\noindent - $D_n := \bigcup_{p\in\omega} D_{n,p}$ soit dense dans $D_{g_n}\cap Z\cap 
g_n^{-1}(T)$.\smallskip

\noindent - $\delta (D_{n,p})$, $\delta (g_n[D_{n,p}])\leq 2^{-\Delta (n,p)}$.

\vfill\eject

 Soit $(x^n_p)_p$ une suite dense de $D_{g_{n}}\cap Z\cap g_n^{-1}(T)$. On choisit un 
ouvert-ferm\'e $D_{n,p}$ de 
$$D_{g_{n}}\cap Z\cap g_n^{-1}(T)$$ 
contenant $x^n_p$, de 
diam\`etre au plus $2^{-\Delta (n,p)}$, dont l'image par $g_n$ soit aussi de diam\`etre au 
plus $2^{-\Delta (n,p)}$. Il est clair que $D_{n,p}$ convient.\bigskip

\noindent $\bullet$ Il reste \`a poser $h_{n,p} := g_{n}\lceil D_{n,p}$. On a ${U_{n,p,U,V} = 
U_{n,U,V}\cap D_{n,p}}$. Enfin, on a vu qu'il y a une infinit\'e d'entiers 
$n$ tels que ${\mbox{Gr}(g_n)\cap\overline{M}\cap (U\times V)\not=\emptyset}$ si ${M\cap 
(U\times V)\not=\emptyset}$, et on a un ouvert non vide $U_{n,U,V}$ de $D^\emptyset_n$. 
Par densit\'e, on trouve $p$ tel que 
${U_{n,U,V}\cap D_{n,p}\not=\emptyset}$. Par 
suite, on a une 
infinit\'e de couples $(n,p)$ tels que ${\mbox{Gr}(h_{n,p})\cap \overline {M}\cap (U\times V)
\not=\emptyset}$.
 Il reste \`a poser $u := \mbox{Id} : Z\rightarrow X$ et $v := \mbox{Id} : 
T\rightarrow Y$.$\hfill\square$\bigskip

 A partir d'un syst\`eme r\'educteur, on n'obtient pas tout de suite une situation 
d'arriv\'ee, mais seulement une situation g\'en\'erale dans un premier temps :
 
\begin{thm} Soit $(Z,T,(k_{q,p})_{(q,p)\in\omega^2},N)$ 
un syst\`eme r\'educteur. Alors il existe $\Phi : \omega^2\rightarrow \omega^2$ injective et des ouverts-ferm\'es $D_{m,p}\subseteq 
D_{k_{\Phi(m,p)}}$ tels que si $g_{m,p} := k_{\Phi(m,p)}\lceil D_{m,p}$,\smallskip

\noindent (a)~$(Z,T,(g_{m,p})_{(m,p)\in\omega^2})$ soit une situation g\'en\'erale.\smallskip

\noindent (b)~Pour $x$ dans $G(g)$ et $y$ dans $\overline{g[x]}\setminus g[x]$, $(x,y)$ est dans 
$N$.\smallskip
 
 Si de plus $(Z,T,(k_{q,p})_{(q,p)\in\omega^2})$ est une 
situation d'arriv\'ee, on peut avoir $G(g)\subseteq G(k)$ et 
$(Z,T,(g_{m,p})_{(m,p)\in\omega^2})$ d'arriv\'ee.\end{thm}

\noindent\bf D\'emonstration.\rm\ Soit 
$$\nu : \left\{\!\!
\begin{array}{ll} 
\omega^{<\omega}\!\!\!\! 
& \rightarrow \omega\cr 
s 
& \mapsto \left\{\!\!\!\!\!\!\!\!
\begin{array}{ll} 
& 0~\mbox{si}~s=\emptyset\mbox{,}\cr 
& \Sigma_{i<\vert s\vert }~(s(i)+1)~\mbox{sinon.}
\end{array}\right.
\end{array}\right.$$
On choisit une distance $d'\leq 1$ compl\`ete sur $N$. On note, pour $x\in Z$, $d'_x$ la 
distance compl\`ete sur $N^x$ d\'efinie par la formule $d'_x(y,t) := d'((x,y),(x,t))$. Pour 
exprimer l'absence de points isol\'es dans $g[x]$, il est plus commode d'indexer les fonctions 
par $\omega^{<\omega}\times \omega$ que par $\omega^2$. 
Soient donc ${e : \omega\rightarrow \omega^{<\omega}}$ bijective, 
$\psi := e\times \mbox{Id}_\omega$ et $\theta := e^{-1}\times \mbox{Id}_\omega$. 
On va en fait construire une injection $\phi :\omega^{<\omega}
\times \omega\rightarrow \omega^2$, et on posera 
ensuite $\Phi := \phi~\circ~\psi$. On va \'egalement construire, par r\'ecurrence sur 
$\vert t\vert $, o\`u $t\in \omega^{<\omega}$, \bigskip

\noindent - Des ouverts-ferm\'es $D_{t,p}$ de $D_{k_{\phi(t,p)}}$.\smallskip
 
\noindent - Des $G_\delta$ denses ${\cal G}_t$ de $D_t := \bigcup_{p\in\omega} 
D_{t,p}$.\smallskip

\noindent - Des ouverts \`a coupes verticales ouvertes-ferm\'ees $\omega_t$ de $Z\times T$.\smallskip
 
\noindent - Des ouverts $G_t$ de $N$.
 
\vfill\eject
 
 On demande \`a ces objets de v\'erifier\smallskip
$$\begin{array}{ll}
& (1)~\overline{D_t} = Z \cr 
& (2)~D_{t,p}\cap D_{t,q}=\emptyset~\mbox{si}~p\not= q\cr 
& (3)~\forall~x\in D_{t^\frown n}\cap {\cal G}_t~~d(k_t(x),k_{t^\frown n}(x))\leq 
2^{1-\nu (t^\frown n)}\mbox{, o\`u}~k_t~\mbox{est~le~recollement~des~fonctions}\cr 
& k_{\phi(t,p)}\lceil D_{t,p}\mbox{, pour}~p~\mbox{entier}\cr 
& (4)~Gr (k_t)\subseteq \omega_t\cap \overline{G_t}\cr 
& (5)~\omega_{t^\frown n}\subseteq \omega_t~\mbox{et}~\omega_{t^\frown n}\cap 
\omega_{t^\frown m}=\emptyset~\mbox{si}~n\not= m\cr 
& (6)~\forall~x\in Z~~\delta (\omega_t^x)\leq 2^{-\nu (t)}\cr 
& (7)~\forall~x\in Z~~\overline{G_{t^\frown n}^x}\cap N^x\subseteq G_t^x \subseteq 
\omega_t^x\mbox{, et}~{\cal G}_t\subseteq \Pi_{Z}[G_t]\cr 
& (8)~\forall~x\in Z~~\delta'_x (G_t^x)\leq 2^{-\vert t\vert }
\end{array}$$
$\bullet$ Admettons ceci r\'ealis\'e. On pose $D_{m,p} := D_{\psi (m,p)}$, $G(g) := 
\bigcap_{t\in \omega^{<\omega}} {\cal G}_t$. Les conditions (a) et (b) 
d'une situation g\'en\'erale sont clairement r\'ealis\'ees. On a $D_{g_{m,p}} = D_{m,p} = 
D_{\psi (m,p)}$, donc ${\bigcup_{p\in\omega} D_{g_{m,p}} = \bigcup_{p\in\omega} 
D_{\psi (m,p)} =\bigcup_{p\in\omega} D_{e(m),p} = D_{e(m)}}$. L'ouvert 
$\bigcup_{p\in\omega} D_{g_{m,p}}$ est donc  
dense dans $Z$, par (1). De plus, la r\'eunion est disjointe, par (2). D'o\`u la 
condition (c) d'une situation g\'en\'erale. On a $G(g) = \bigcap_{t\in \omega^{<\omega}} 
{\cal G}_t$, donc $G(g)$ est $G_\delta$ dense de $Z$, et si $x\in 
G(g)$ et $m\in\omega$, 
$${x\in {\cal G}_{e(m)}\subseteq D_{e(m)} = 
\bigcup_{p\in\omega} D_{g_{m,p}} = D_{g_m}}.$$ 
Donc $G(g)\subseteq\bigcap_{m\in\omega} D_{g_m}$. Par (3) et la remarque qui suit, on a la condition (d) 
d'une situation g\'en\'erale.\bigskip

 Notons que ${g[x] = \{k_t(x)~/~t\in 
\omega^{<\omega}\}}$. En effet, si $t\in \omega^{<\omega}$, $\exists~p\in\omega$ tel que 
$${k_t(x) = k_{\phi (t,p)}(x) = k_{\Phi (\theta (t,p))}(x) = g_{\theta (t,p)}(x) = g_{e^{-1}(t)}(x)}.$$ 
Inversement, si $m\in\omega$, on trouve $p\in\omega$ tel que 
${g_m(x) = g_{m,p}(x) = k_{\Phi (m,p)}(x) = k_{e(m)}(x)}$.\bigskip

 L'argument qui suit est semblable \`a celui utilis\'e dans la preuve du th\'eor\`eme 
d'Hurewicz dans [SR]. On pose $M_k := \{ k_t(x)~/~\vert t\vert  \leq k\}$. Alors 
$M_k$ est ferm\'e dans $T$, par r\'ecurrence sur $k$ : si $M_k^\varepsilon := \{ y\in 
T~/~d(y,M_k)\leq \varepsilon\}$, on a $M_{k+1} = \bigcap_{\varepsilon > 0} [M_k^\varepsilon \cup 
(M_{k+1}\setminus M_k^\varepsilon)]$, et $M_{k+1}\setminus M_k^\varepsilon$ est fini, par (3).\bigskip

 Si $k\in\omega$ et $y\in \overline{g[x]}\setminus g[x]$, ${y\notin M_k}$. Donc il 
existe ${\varepsilon > 0}$ tel que ${y\notin M_k^\varepsilon}$. Si ${t\in \omega^k}$ et 
$(x,w)$ est dans $\omega_t$, ${d(w,M_k)\leq d(w,k_t(x))\leq {\it\Delta}_x (\omega_t^x)\leq 
2^{-\nu (t)} \leq \varepsilon}$ d\`es que $\nu (t)\geq k_0$, donc $\omega_t^x \subseteq M_k^
\varepsilon$, sauf pour un nombre fini de $t$ dans $\omega^k$. Donc si on 
d\'efinit ${\cal H} := \{t\in \omega^k~/~\omega_t^x \not\subseteq 
M_k^\varepsilon\}$, on a la suite d'inclusions suivante : $${g[x]\! =\!  
M_k \! \cup\!  \{ k_t(x)~/~\vert t\vert \! >\! k\}\!  \subseteq\!  M_k\!  \cup 
\bigcup_{\vert w\vert\! >\! k} \omega_w^x \! \subseteq\!  M_k\!  \cup \bigcup_{\vert t\vert \! 
 =\!  k} \omega_t^x \! \subseteq \! M_k^\varepsilon \! \cup \bigcup_{t\in {\cal H}} 
\omega_t^x}.$$ 
D'o\`u $\overline{g[x]} \subseteq M_k^\varepsilon \cup \bigcup_{t\in {\cal H}} 
\omega_t^x \subseteq M_k^\varepsilon \cup \bigcup_{\vert t\vert  = k} \omega_t^x$. Donc on trouve 
une unique suite $\sigma$ dans $\omega^\omega$ telle que 
$y\in \bigcap_{t\prec \sigma} \omega_t^x$. $(\overline{G_t^x}
\cap N^x)_{t\prec \sigma}$ est une suite d\'ecroissante de ferm\'es non vides 
dont les diam\`etres tendent vers 0 de $(N^x, d'_x)$, par (7) et (8), donc converge vers 
$\xi\in N^x$, et $\{\xi\} = \bigcap_{t\prec\sigma} G_t^x$. 
D'o\`u $\xi \in \bigcap_{t\prec\sigma} \omega_t^x = \{y\}$ et $(x,y)\in N$.

\vfill\eject

\noindent $\bullet$ Montrons donc que la construction est possible. Soit $(Z_n)$ une base de la 
topologie de $Z$ form\'ee d'ouverts-ferm\'es non vides. Comme $Z$ est polonais de dimension 0, il 
peut \^etre vu comme un ferm\'e de $\omega^\omega$ ; on peut donc supposer que $Z$ est muni 
d'une distance compl\`ete telle que $d(Z_n,\check Z_n)>0$. On pose $\omega_\emptyset := 
Z\times T$, $G_\emptyset := N$. On construit $\phi (\emptyset,p)$ et $D_{\emptyset,p}$ 
par r\'ecurrence sur $p$, en exigeant que $\bigcup_{q<p} D_{\emptyset,q}
\not= Z$. Admettons avoir construit $\phi (\emptyset,q)$ et 
$D_{\emptyset,q}$ pour $q<p$. On peut d\'efinir, pour $q\leq p$, 
$$n(q) := \mbox{min}\{n\in\omega~/~Z_n\cap\bigcup_{r<q} D_{\emptyset ,r}=\emptyset~\mbox{et}~\forall~r<q~~n>n(r)\}.$$ 
Comme $\Pi_Z[N]$ est dense dans $Z$, $N\cap (Z_{n(p)}\times T)\not=
\emptyset$, donc par la condition (e) d'un syst\`eme r\'educteur, on peut trouver $(m,r)\in
\omega^2\setminus\{\phi (\emptyset ,q)~/~q<p\}$ tel que 
$\delta (D_{k_{m,r}}) < d(Z_{n(p)},\check Z_{n(p)})$ et $\mbox{Gr}(k_{m,r})\cap \overline{N}\cap 
(Z_{n(p)}\times T)\not=\emptyset$. On a donc $D_{k_{m,r}}\subseteq Z_{n(p)}$, et on pose 
$\phi (\emptyset, p):=(m,r)$ ; on choisit $D_{\emptyset ,p}\subseteq U_{m,r,Z_{n(p)},T}$ tel 
que $\bigcup_{q\leq p} D_{\emptyset,q}\not= Z$. De sorte que la condition 
(1) est r\'ealis\'ee. Ceci termine la construction pour $\vert t\vert  =0$.\bigskip
 
\noindent $\bullet$ On effectue maintenant une sous-construction : on construit, par r\'ecurrence sur $n$,\bigskip

\noindent - Les ouverts $\omega_{t^\frown n}$.\smallskip

\noindent - Une suite d\'ecroissante $(E_n)$ de $G_\delta$ denses de $D_t$. \smallskip

\noindent - Des fonctions continues $f_n : E_n\rightarrow T$.\smallskip

\noindent - Des ouverts $V_n$ de $Z\times T$.\bigskip
 
\noindent On demande \`a ces objets de v\'erifier 
$$\begin{array}{ll}
& (i)~~~\mbox{Gr}(f_n)\subseteq \omega_{t^\frown n}\cap G_t\cr 
& (ii)~~\omega_{t^\frown n}\subseteq\omega_t~\mbox{et}~\omega_{t^\frown n}\cap
\omega_{t^\frown m} = \emptyset~\mbox{si}~n\not= m\cr 
& (iii)~\forall~x\in E_n~~d(f_n(x),k_t(x))\leq 2^{-\nu (t^\frown n)} \cr 
& (iv)~~\forall~x\in Z~~\delta (\omega_{t^\frown n}^x)\leq 2^{-\nu (t^\frown n)} \cr 
& (v)~~~\mbox{Gr}(k_t\lceil E_n)\subseteq V_n\subseteq (Z\times T)\setminus (\bigcup_{m\leq n} 
\omega_{t^\frown m})
\end{array}$$
$\bullet$ Montrons qu'une telle construction est possible. Si on a construit ces objets pour 
$q<n$, soit $x\in D_t$. Par continuit\'e de $k_t$, on peut trouver un 
voisinage ouvert ${\cal W}_x$ de $x$ inclus dans $D_t$ et un voisinage 
ouvert-ferm\'e ${\cal V}_x$ de diam\`etre au plus $2^{-\nu (t^\frown n)}$ de $k_t(x)$ 
tels que ${\cal W}_x\subseteq k_t^{-1}
({\cal V}_x)$. De sorte que si $(z,t)\in {\cal W}_x\times {\cal V}_x$, $d(t,k_t(z))\leq 
2^{-\nu (t^\frown n)}$. On a, par la propri\'et\'e de Lindel\"of, 
$$\bigcup_{x\in D_t}({\cal W}_x\times {\cal V}_x) = \bigcup_{m\in\omega} ({\cal W}_m\times {\cal V}_m).$$  
On r\'eduit la suite $({\cal W}_m)$ en $({\cal W}'_m)$, puis on pose 
${\cal U}_n := \bigcup_{m\in\omega} ({\cal W}'_m\times {\cal V}_m)$. Soit 
${\cal O}_n$ un ouvert dense de 
$D_t$ contenant $E_{n-1}$ tel que $\mbox{Gr}(k_t) \cap V_{n-1} = \mbox{Gr}(k_t\lceil 
{\cal O}_n)$. Par le th\'eor\`eme de Jankov-von Neumann, on peut trouver $\tilde f_n$ 
Baire-mesurable uniformisant ${\cal U}_n\cap G_t$ sur sa projection $\Pi$. 
$\Pi$ est dense dans $Z$. En effet, soient $U$ un ouvert-ferm\'e non vide de $Z$ et 
$x\in {\cal O}_n\cap U$. Comme $(x,k_t(x))\in {\cal U}_{n}\cap
\overline{G_t}$, on peut trouver un voisinage ouvert-ferm\'e $W$ de $x$ tel que 
$W\times k_t[W]\subseteq {\cal U}_{n}\cap (U\times T)$. Soit alors $(z,y)
\in (W\times k_t[W])\cap G_t$. On a que $z\in U\cap\Pi\not=\emptyset$. Par la 
condition (d) d'un syst\`eme r\'educteur, $\Pi$ est comaigre dans $Z$, donc contient un 
$G_\delta$ dense $W_n$ de $Z$. Alors $\tilde f_n\lceil W_n$ est Baire-mesurable, donc 
on peut trouver un $G_\delta$ dense $F_n$ de $W_n\cap D_t$ tel que $\tilde f_n\lceil F_n$ soit continue.

\vfill\eject

 On peut poser $E_n := F_n\cap E_{n-1}$ et $f_n := \tilde f_n\lceil E_n$. Si $x\in E_n$, 
$(x,f_n(x))\in {\cal U}_n$, donc $d(f_n(x), k_t(x))\leq 2^{-\nu (t^\frown n)}$. Les 
graphes de $f_n$ et $k_t\lceil E_n$ sont des ferm\'es disjoints de $E_n\times T$, donc 
on peut trouver un ouvert-ferm\'e $\theta$ de $E_n\times T$ tel que $\mbox{Gr}(f_n)\subseteq 
\theta\subseteq (E_n\times T)\setminus \mbox{Gr}(k_t\lceil E_n)$.\bigskip

 On peut trouver des ouverts disjoints $\cal T$ et $\cal W$ de $Z\times T$ tels que 
$\theta = (E_n\times T)\cap {\cal T}$ et 
$$(E_n\times T)\setminus \theta = (E_n\times T)\cap {\cal W}\mbox{,}$$
par la propri\'et\'e de r\'eduction des ouverts. Soit $(T_m)$ une base de la topologie de $T$ stable 
par intersections finies et form\'ee d'ouverts-ferm\'es v\'erifiant $d(T_m,\check T_m)>0$. 
On raisonne alors comme pr\'ec\'edemment : pour 
$x$ dans $E_n$, on trouve un voisinage ouvert-ferm\'e de base ${\cal V}'_x$ de $f_n(x)$ de 
diam\`etre au plus $2^{-\nu (t^\frown n)}$ et un voisinage ouvert-ferm\'e ${\cal Y}_x$ de $x$ 
tels que ${\cal Y}_x\cap E_n\subseteq f_n^{-1}({\cal V}'_x)$ et ${\cal Y}_x\times 
{\cal V}'_x\subseteq {\cal T}\cap V_{n-1}\cap \omega_t$. Comme avant, on applique la 
propri\'et\'e de Lindel\"of, ce qui fournit ${\cal Y}_m$ et ${\cal V}'_m$, et on r\'eduit 
la suite $({\cal Y}_m)$ en $({\cal Y}'_m)$. On pose $\omega_{t^\frown n} := 
\bigcup_{m\in\omega} {\cal Y}'_m\times {\cal V}'_m$ et $V_n := V_{n-1}\cap {\cal W}$. 
Les conditions (i) \`a (v) sont clairement satisfaites. On a donc les conditions (5) 
et (6) de la construction principale.\bigskip

\noindent $\bullet$ On proc\`ede encore comme avant pour d\'efinir $G_{t^\frown n}$. Pour $x$ dans 
$E_n$, on trouve un voisinage ouvert-ferm\'e de base ${\cal D}_x$ de $f_n(x)$ et un 
voisinage ouvert-ferm\'e ${\cal C}_x$ de $x$ tels que
$${\delta' ([{\cal C}_x\times {\cal D}_x]\cap N)\leq 2^{-\vert t\vert -1}}\mbox{,}$$
${{\cal C}_x\cap E_n\subseteq f_n^{-1}({\cal D}_x)}$ et 
${{\cal C}_x\times {\cal D}_x\subseteq\omega_{t^\frown n}\cap G\cap (D_t
\times T)}$, o\`u $G$ est ouvert de $Z\times T$ tel que 
$${G\cap N = G_t}.$$ 
On applique la propri\'et\'e de Lindel\"of, ce qui fournit ${\cal C}_m$ et 
${\cal D}_m$, et on r\'eduit la suite ${({\cal C}_m)}$ en ${({\cal C}'_m)}$. On pose 
${G_{t^\frown n} := N\cap (\bigcup_{m\in\omega} {\cal C}'_m\times {\cal D}_m)}$ et 
${{\cal G}_t := \bigcap_{n\in\omega} E_n}$, et les conditions (7) et (8) sont satisfaites. 
On a ${G_{t^\frown n}\cap\omega_{t^\frown n} = N\cap (\bigcup_{m\in\omega} 
{\cal C}'_m\times {\cal D}_m)\cap (\bigcup_{m\in\omega} {\cal Y}'_m\times 
{\cal V}'_m) = \bigcup_{l\in\omega} N\cap (A_l\times B_l)}$, avec par d\'efinition  
${A_l := {\cal C}'_{e_2^{-1}(l)_0}\cap {\cal Y}'_{e_2^{-1}(l)_1}}$ et 
${B_l := {\cal D}_{e_2^{-1}(l)_0}\cap {\cal V}'_{e_2^{-1}(l)_1}}$, o\`u ${e_2 :\omega^2
\rightarrow\omega}$ est bijective. De plus, on a 
${A_l\cap A_{l'} = \emptyset}$ si $l\not= l'$ et ${\Pi_Z[N\cap (A_l\times B_l)]}$ est 
comaigre dans $A_l$ car ${\mbox{Gr}(f_n)\subseteq G_{t^\frown n}\cap \omega_{t^\frown n}}$. 
On proc\`ede alors comme pour $t=\emptyset$ pour terminer la construction, en travaillant dans 
${A_l\times B_l}$. On construit, par r\'ecurrence sur $p$, ${\phi (t^\frown n,e_2(l,p))}$ et 
${D_{t^\frown n,e_2(l,p)}}$, en exi-geant que ${\bigcup_{q<p} 
D_{t^\frown n,e_2(l,q)}\not= A_l}$. Soit $(Z_{n_k})_{k}$ la base de 
$A_l$ form\'ee des $Z_n\subseteq A_l$. On pose
$$k(q) := \mbox{min}\{k\in\omega~/~Z_{n_k}\cap \bigcup_{r<q} 
D_{t^\frown n,e_2(l,r)}=\emptyset~\mbox{et}~\forall~r<q~~k>k(r)\}.$$
On choisit $(m,r)\in \omega^2\setminus \{\phi (t^\frown n,e_2(l,q))~/~
q<p\}$ tel que $\delta (D_{k_{m,r}}) < d(Z_{n_{k(p)}},\check Z_{n_{k(p)}})$, 
$${\delta (\mbox{Im}(k_{m,r})) < d(B_l,\check B_l)}\mbox{,}$$ 
et $\mbox{Gr}(k_{m,r})\cap \overline{N}\cap (Z_{n_{k(p)}}\times B_l)\! \not=\! \emptyset$. On pose 
${\phi (t^\frown n,e_2(l,p))\! :=\!  (m,r)}$ et on choisit $D_{t^\frown n,e_2(l,p)}$ dans 
$U_{m,r,Z_{n_{k(p)}},B_l}$ tel que $\overline {\bigcup_{q\leq p} D_{t^\frown n,e_2(l,q)}}\not= A_l$. Ceci termine la construction car 
$$\mbox{Gr}(k_{m,r})\subseteq Z_{n_{k(p)}}\times B_l\subseteq \omega_{t^\frown n}\mbox{,}$$ 
ce qui assure l'injectivit\'e de $\phi$.\bigskip

\noindent $\bullet$ Si de plus $(Z,T,(k_{q,p})_{(q,p)\in\omega^2})$ est d'arriv\'ee, quitte \`a 
remplacer $G(g)$ par $G(g)\cap G(k)$, on peut avoir $G(g) \subseteq G(k)$. 
L'injectivit\'e de $\Phi$ fait que la situation g\'en\'erale $(Z,T,(g_{m,p})_{(m,p)\in
\omega^2})$ est en fait d'arriv\'ee.$\hfill\square$

\vfill\eject

 On doit maintenant obtenir une situation d'arriv\'ee \`a partir d'une situation g\'en\'erale :

\begin{thm} Soit $(Z,T,(l_{r,p})_{(r,p)\in\omega^2})$ une situation g\'en\'erale. Alors il existe une 
situation d'arriv\'ee $(Z,T,(k_{q,p})_{(q,p)\in\omega^2})$ telle que $G(k)\subseteq G(l)$ et pour $x$ dans 
$G(k)$, on ait ${k[x]\subseteq l[x]}$.\end{thm}

\noindent\bf D\'emonstration.\rm\ Pour exprimer \`a la fois 
l'absence de points isol\'es dans $k[x]$, la densit\'e de $D_{k_q}$ dans $Z$ et la 
condition (c) d'une situation d'arriv\'ee, il est plus commode d'indexer les fonctions 
par $\bigcup_{n\in\omega} \omega^n\times\omega^{n+1}$ que par $\omega^2$. Soient donc 
$\psi\! :\!\mbox{Im}(M)\!\rightarrow\!\omega$ bijective croissante, et ${\phi := 
\psi\circ M}$, o\`u 
$$M : \left\{\!\!
\begin{array}{ll}
\bigcup_{n\in\omega} \omega^n\times\omega^{n+1}\!\!\!\! 
& \rightarrow \omega \cr 
(s,t) 
& \mapsto  q_0^{t(0)+1}q_1^{s(0)+1}...q_{2n-2}^{t(n-1)+1}q_{2n-1}^{s(n-1)+1}q_{2n}^{t(n)+1}
\end{array}\right.$$
($(q_n)$ est la suite des nombres premiers). $\phi$ est donc une bijection de 
$\bigcup_{n\in\omega} \omega^n\times\omega^{n+1}$ sur $\omega$ v\'erifiant 
$\phi (s,t)<\phi (s^\frown p,t^\frown m)$ et $\phi (s,t) < 
\phi (s,t\lceil \vert s\vert ^\frown [t(\vert s\vert )+1])$ pour $(s,t)$ dans $\bigcup_{n\in\omega} 
\omega^n\times\omega^{n+1}$ et $(p,m)$ dans $\omega^2$. Soient $(Z_n)$ une 
base de la topologie de $Z$, et $(O_r)$ une suite d'ouverts denses de $Z$ tels que 
$G(l) = \bigcap_{r\in\omega} O_r$. En supposant 
$$\forall~n\in\omega~~\forall~s,t\in\omega^n~~\forall~m\in 
\omega~~D_{k_{s\lceil (\vert s\vert -1),t}}\! \setminus  
(\bigcup_{i<m} D_{k_{s,t^\frown i}})\! \not=\!  \emptyset$$ 
avec la convention $D_{k_{\emptyset\lceil -1,\emptyset}} := Z$, on d\'efinit 
${r :\bigcup_{n\geq -1} \omega^n\times\omega^{n+1}\rightarrow \omega\cup\{ -1\}}$  en posant 
$$\begin{array}{ll}
r(\emptyset\lceil -1,\emptyset) & \!\!\!\! := -1\mbox{,}\cr 
r(s,t^\frown m) & \!\!\!\! :=\mbox{min}\left\{ r\in\omega~/~\left\{\!\!\!\!\!\!
\begin{array}{ll} 
& Z_r\cap D_{k_{s\lceil (\vert s\vert -1),t}}\! \setminus (\bigcup_{i< m} D_{k_{s,t^\frown i}})\! \not=\!\emptyset\cr 
& \mbox{et}\cr 
& r\!>\! \mbox{max}(r(s\lceil (\vert s\vert -1),t),\mbox{max}_{i< m}~r(s,t^\frown i))
\end{array}
\right.\right\}.
\end{array}$$
On a donc $r(s,t^\frown (n+1))> r(s,t^\frown n)>r(s\lceil (\vert s\vert -1),t)$. On 
construit des fonctions $k_{s,t^\frown m}$, pour $(s,t,m)\in (\bigcup_{n\in\omega} 
(\omega^n)^2)\times\omega$, par r\'ecurrence sur $\phi (s,t^\frown m)$, en demandant 
$$\begin{array}{ll}
& (1)~D_{k_{s,t^\frown m}}\subseteq Z_{r(s,t^\frown m)}\cap D_{k_{s\lceil (\vert s\vert -1),t}}\setminus 
(\bigcup_{i< m} D_{k_{s,t^\frown i}})\cap O_{\vert s\vert }\cr 
& (2)~\bigcup_{i< m} D_{k_{s,t^\frown i}}\not= D_{k_{s\lceil (\vert s\vert -1),t}}~\mbox{et}~
D_{k_{s\lceil (\vert s\vert -1),t}}\not=\emptyset\cr 
& (3)~\delta (D_{k_{s,t^\frown m}}),~\delta (\mbox{Im}(k_{s,t^\frown m})) \leq 2^{-\phi 
(s,t^\frown m)}\cr 
& (4)~k_{s,t^\frown m}~\mbox{est~la~restriction~d}\hbox {'}\mbox{une~des~fonctions}~l_{r,p}
~\mbox{\`a~un~ouvert}\hbox {-}\mbox{ferm\'e~de}~D_{l_{r,p}}\cr 
& (5)~\forall~n,m\in\omega~~\forall~s,t\in\omega^n~~\forall~x\in D_{k_{s,t^\frown m}}~~d(k_{s,t^\frown m}(x),k_{s\lceil (\vert s\vert -1),t}(x))<2^{-s(\vert s\vert -1)}\cr 
& (6)~\forall~h\in\omega\setminus\{ 0\}~~\forall~u\in (\bigcup_{p\leq \phi(s,t^\frown m)} 
\{\phi^{-1}(p)\})^{2h}\cr\cr 
& ~~~\pmatrix { 
& \exists~U\in\borone,~\emptyset\!\not=\! U\subseteq Z\cr 
& \mbox{et}\cr 
& \!\!\!\!\!\!\forall~x\in U~~k_{u(0)}^{-1}k_{u(1)}^{}...k_{u(2k-2)}^{-1}k_{u(2k-1)}^{}(x)\! =\!x}~\Rightarrow ~
\exists~i\! <\! 2h-1~~u(i)\! =\! u(i+1)\cr\cr 
& (7)~\mbox{Il~n}\hbox {'}\mbox{y~a~qu}\hbox {'}\mbox{un~nombre~fini~de~compositions~du~type~de~la~condition~(6)~sans~termes}\cr 
& \mbox{cons\'ecutifs~identiques~ayant~un~domaine~de~d\'efinition~
non~vide.}
\end{array}$$

\vfill\eject

\noindent $\bullet$ Admettons cette construction r\'ealis\'ee. Soient 
$e : \omega \rightarrow\omega^{<\omega}$ et, pour $n>0$, $e_n : \omega^n \rightarrow \omega$ bijectives. On pose $k_{q,p} := k_{e(q),e_{\vert e(q)\vert+1}^{-1}(p)}$. Les applications 
suivantes sont r\'eciproques l'une de l'autre : 
$$\left\{\!\!
\begin{array}{ll}
\omega^2 
& \!\!\!\!\rightarrow \bigcup_{n\in\omega} \omega^n\times\omega^{n+1} \cr 
(q,p) 
& \!\!\!\!\mapsto (e(q),e_{\vert e(q)\vert+1}^{-1}(p))
\end{array}
\right.\mbox{,}~~~\left\{\!\!
\begin{array}{ll}
\bigcup_{n\in\omega}\omega^n\times\omega^{n+1} 
& \!\!\!\!\rightarrow\omega^2 \cr 
(s,t^\frown m) 
& \!\!\!\!\mapsto (e^{-1}(s),e_{\vert s\vert +1}(t^\frown m))
\end{array}
\right..$$  
Les conditions (a) et (b) d'une situation g\'en\'erale et les conditions (b) et (c) d'une situation d'arriv\'ee seront alors clairement r\'ealis\'ees, par (3), (4) et (6).  $\bigcup_{m\in\omega} D_{k_{s,t^\frown m}}$ 
est dense dans $D_{k_{s\lceil (\vert s\vert \! -\! 1),t}}$, sinon on peut trouver $r\! >\! 
r(s\lceil (\vert s\vert \! -\! 1),t)$ tel que ${Z_r\! \subseteq \! 
D_{k_{s\lceil (\vert s\vert \! -\! 1),t}}\! \setminus\! \overline{\bigcup_{m\in\omega} D_{k_{s,t^
\frown m}}}}$. Comme la suite $(r(s,t^\frown m))_m$ cro\^\i t strictement vers l'infini, on trouve un plus petit 
$m$ tel que $r(s,t^\frown m)>r$.\bigskip

 Si $m>0$, $r(s,t^\frown (m-1)) \leq r$ ; si $r(s,t^\frown (m-1))<r$, on a $r(s,t^\frown m) \leq r$, ce 
qui est absurde. Si $r = r(s,t^\frown (m-1))$, 
${D_{k_{s,t^\frown (m-1)}}\subseteq Z_r}$, ce qui contredit la disjonction de $Z_r$ et 
$D_{k_{s,t^\frown (m-1)}}$. Si $m=0$, on a $r(s,t^\frown 0)>r>r(s\lceil (\vert s\vert -1),t)$, ce qui contredit la 
d\'efinition de $r(s,t^\frown 0)$. \bigskip 

 Il suffit donc d'assurer (1) et (2) pour avoir la disjonction de $D_{k_{s,t^\frown m}}$ et 
$D_{k_{s,t^\frown n}}$ pour $n\not= m$, donc de $D_{k_{s,t^\frown m}}$ et 
$D_{k_{s,v^\frown n}}$ pour $(t,m)\not= (v,n)$, et la densit\'e de ${D_s := 
\bigcup_{(t,m)\in\omega^{\vert s\vert }\times\omega} D_{k_{s,t^\frown m}}}$ dans 
$\bigcup_{t\in\omega^{\vert s\vert }} D_{k_{s\lceil (\vert s\vert -1),t}}$ ; on a donc la condition (c) d'une situation g\'en\'erale. Enfin, on pose 
$$G(k) := \bigcap_{q\in\omega} D_{k_q}$$ 
et la condition (5) entra\^\i ne la condition (d) d'une situation g\'en\'erale car si $q\in\omega$ et 
$$x\in G(k) = \bigcap_{s\in\omega^{<\omega}} D_s\mbox{,}$$
$x\in D_{e(q)}\cap D_{e(q)^\frown p}$ pour tout entier $p$, donc il existe $t\in\omega^{\vert e(q)\vert+1}$ et 
$m\in\omega$ tels que 
$${x\in D_{e(q)^\frown p,t^\frown m}\subseteq D_{e(q),t}}\mbox{,}$$ 
et ${d(k_q(x),k_l(x))<2^{-p}}$, o\`u $e(l) = e(q)^\frown p$.\bigskip

\noindent $\bullet$ Montrons donc que cette construction est possible. Soit $(r,p)\in\omega^2$ 
telle que $D_{l_{r,p}}\cap Z_0\not=\emptyset$. On choisit 
un ouvert-ferm\'e non vide $D_{k_{\phi^{-1}(0)}}$ strictement inclus dans 
$D_{l_{r,p}}\cap Z_0\cap O_{\vert\phi_0^{-1}(0)\vert }$ et on pose $k_{\phi^{-1}(0)} := 
l_{r,p}\lceil D_{k_{\phi^{-1}(0)}}$. Admettons avoir 
construit $(k_{s,t^\frown m})_{\phi (s,t^\frown m)\leq n}$ v\'erifiant (1)-(7), ce qui est fait pour $n=0$. Posons 
${\phi^{-1}(n+1) := (s,t^\frown m)}$. On a d\'ej\`a construit $D_{k_{\phi^{-1}(q_1)}}$, ..., 
$D_{k_{\phi^{-1}(q_m)}}$ dans $D_{k_{s\lceil (\vert s\vert -1),t}}$, de sorte que 
${\phi^{-1}(q_i) = (s,t^\frown (i-1))}$ pour $1\leq i\leq m$, par construction de $\phi$.\bigskip

 Soit $U$ un ouvert-ferm\'e non vide de $D_{k_{s\lceil (\vert s\vert -1),t}}$ tel 
que $\delta (U)$, $\delta (k_{s\lceil (\vert s\vert -1),t}[U])\leq 2^{-n-1}$, 
$${U\cup \bigcup_{i< m} D_{k_{s,t^\frown i}} \not= D_{k_{s\lceil (\vert s\vert -1),t}}}\mbox{,}$$ 
et tel que $U\times k_{s\lceil (\vert s\vert -1),t}[U]$ soit inclus dans 
$${\left\{(x,y)\in [Z_{r(s,t^\frown m)}\cap 
D_{k_{s\lceil (\vert s\vert -1),t}}\setminus 
(\bigcup_{i< m} D_{k_{s,t^\frown i}})]\times T~/~d(y,k_{s\lceil (\vert s\vert -1),t}(x))<
2^{-s(\vert s\vert -1)}\right\}}.$$

\vfill\eject
 
 Ce qui suit est \`a rapprocher du lemme 2.11 de [Le3]. Par (7), on note $\{ H_1, ...,H_p\}$ l'ensemble fini des compositions du type de la condition (6) sans termes cons\'ecutifs identiques de domaine 
de d\'efinition non vide d\'efinissables \`a partir de $k_{\phi^{-1}(0)}$, ..., 
$k_{\phi^{-1}(n)}$. Pour $I\subseteq p$, on pose 
$${O_I := \bigcap_{i\in I} D_{H_{i+1}}\cap\bigcap_{i\in p\setminus I}\check D_{H_{i+1}}}.$$
Alors $(O_I)_{I\subseteq p}$ est une partition en ouverts-ferm\'es de $Z$, donc on trouve $I\subseteq p$ tel que $U\cap O_I \not= \emptyset$. Soit ${O' := \{x\in U\cap O_I~/~\forall~i\in I~~H_{i+1}(x)\not= 
x\}}$. Par (6), $O'$ est un ouvert dense de $U\cap O_I$ ; on peut donc trouver un 
ouvert-ferm\'e non vide $O''$ de $O'$ tel que $\forall~i\in I$, 
${O''\cap H_{i+1}[O''] = \emptyset}$. Soit ${x\in O''\cap G(l)}$, et $(r,p)\in\omega^2$ telle que 
${(x,l_{r,p}(x))\in O''\times k_{s\lceil (\vert s\vert -1),t}[O'']}$, 
$${l_{r,p} (x)\not= k_{\phi^{-1}(j)}H_{m+1}(x)}$$ 
et aussi ${l_{r,p} (x)\not= k_{\phi^{-1}(j)}(x)}$ pour tout $m<p$ et tout ${j\leq n}$. On choisit 
$D_{k_{\phi^{-1}(n+1)}}$ contenant $x$  dans l'ouvert 
${D_{l_{r,p}}\cap O''\cap l_{r,p}^{-1}(k_{s\lceil (\vert s\vert -1),t}[O''])\cap 
O_{\vert\phi_0^{-1}(n+1)\vert }}$ tel que pour tout $m<p$ et pour tout $j\leq n$ 
on ait ${l_{r,p}[D_{k_{\phi^{-1}(n+1)}}] \cap k_{\phi^{-1}(j)}H_{m+1}
[D_{k_{\phi^{-1}(n+1)}}]\! =\! \emptyset}$ et ${l_{r,p}[D_{k_{\phi^{-1}(n+1)}}] \cap 
k_{\phi^{-1}(j)}[D_{k_{\phi^{-1}(n+1)}}] \! =\! \emptyset}$. On pose ${k_{\phi^{-1}(n+1)} 
:= l_{r,p}\lceil D_{k_{\phi^{-1}(n+1)}}}$.\bigskip

  Il est clair que les conditions (1) \`a (5) sont r\'ealis\'ees. Montrons que la condition (6) 
est r\'ealis\'ee au rang $n+1$, en raisonnant par l'absurde. On peut donc trouver $h>0$, 
$u\in (\bigcup_{p\leq n+1} \{\phi^{-1}(p)\})^{2h}$ et un ouvert-ferm\'e non vide $U$ de 
$Z$ tels que pour tout $x$ de $U$ on ait $k_{u(0)}^{-1}k_{u(1)}^{}...k_{u(2h-2)}^{-1}
k_{u(2h-1)}^{}(x) = x$ et $u(i)\not= u(i+1)$ si $i<2h-1$. Par hypoth\`ese de r\'ecurrence, 
on peut trouver $i<2h$ minimal tel que $u(i) = \phi^{-1}(n+1)$. Montrons, en 
raisonnant par l'absurde, qu'un tel $i$ est unique. Si tel n'est pas le cas, on peut 
trouver $j>i+1$ minimal tel que $u(j) = \phi^{-1}(n+1)$. Il y a alors quatre cas.\bigskip

\noindent\bf  Cas 1.~\rm  $i$ et $j$ sont impairs.\bigskip

 Posons $y := k_{u(j)}^{}k_{u(j+1)}^{-1}...k_{u(2h-1)}^{}(x)$. Comme $u(j)=\phi^{-1}(n+1)$, 
$y$ est dans $k_{s\lceil (\vert s\vert -1),t}[O'']$ et on trouve $z$ dans $O''$ tel que 
$y = k_{s\lceil (\vert s\vert -1),t}(z)$. On 
a la double \'egalit\'e 
$${x = k_{u(0)}^{-1}k_{u(1)}^{}...k_{u(j-1)}^{-1}(y) = 
k_{u(0)}^{-1}k_{u(1)}^{}...k_{u(j-1)}^{-1}k_{s\lceil (\vert s\vert -1),t}(z)}\mbox{,}$$ 
donc le domaine d\'efinition de ${k_{u(i+1)}^{-1}k_{u(i+2)}^{}...k_{u(j-1)}^{-1}
k_{s\lceil (\vert s\vert -1),t}^{}}$ est non vide. Si $i\! +\! 1\! =\! j\!-\! 1$ et ${
u(j-1) = s\lceil (\vert s\vert -1),t}$, ${z\in D_{k_{\phi^{-1}(n+1)}}}$ et 
$${y\in k_{\phi^{-1}(n+1)}[D_{k_{\phi^{-1}(n+1)}}]
\cap k_{s\lceil (\vert s\vert -1),t}[D_{k_{\phi^{-1}(n+1)}}] = \emptyset}.$$ 
Si $i+1 \not= j-1$ ou $u(j\! -\! 1)\! \not=\!  s\lceil (\vert s\vert -1),t$, on peut 
trouver $r\! <\! p$ tel que 
$${H_{r+1}\!  =\!  k_{u(i+1)}^{-1}k_{u(i+2)}^{}...k_{u(j\! -\! 1)}^{-1}k_{s\lceil (\vert s\vert -1),t}^{}}\mbox{,}$$ 
comme $z\in O''\subseteq O_I$, ${H_{r+1}(z)\notin O''}$, donc 
${k_{u(i)}(H_{r+1}(z))}$ n'est pas d\'efini, ce qui est absurde.

\vfill\eject

\noindent\bf  Cas 2.~\rm  $i$ et $j$ sont pairs.\bigskip

 Dans la composition appara\^\i t 
$k_{\phi^{-1}(n+1)}^{-1}k_{u(i+1)}^{}...k_{u(j-1)}^{}k_{\phi^{-1}(n+1)}^{-1}$, donc la composition 
$$k_{\phi^{-1}(n+1)}^{}k_{u(j-1)}^{-1}...k_{u(i+1)}^{-1}k_{\phi^{-1}(n+1)}^{}$$ 
a un domaine de d\'efinition non vide. Mais on voit comme avant que c'est impossible.\bigskip

\noindent\bf  Cas 3.~\rm  $i$ est impair et $j$ est pair.\bigskip

 Soit $r<p$ tel que $H_{r+1} = k_{u(i+1)}^{-1}...k_{u(j-1)}^{}$. Alors 
$k_{u(j)}^{-1}...k_{u(2h-1)}^{}(x)\in O''$, donc comme avant 
$k_{u(i+1)}^{-1}...k_{u(2h-1)}^{}(x)\notin O''$, ce qui est absurde.\bigskip

\noindent\bf  Cas 4.~\rm  $i$ est pair et $j$ est impair.\bigskip

 Posons $y := k_{u(j)}^{}...k_{u(2h-1)}^{}(x)$ ; comme dans le cas 1, on trouve $z$ dans $O''$ tel que 
$${y\!  =\!  k_{s\lceil (\vert s\vert -1),t}(z)}.$$ 
Posons $w \! :=\! k_{u(i+2)}^{-1}...k_{u(2h\! -\! 1)}^{}(x)$ ; on a 
$w\!  =\!  k_{u(i+2)}^{-1}...k_{u(j\! -\! 1)}^{-1}k_{s\lceil (\vert s\vert -1),t}^{}(z)$, et 
$$k_{u(i+1)}(w) = k_{s\lceil (\vert s\vert -1),t}^{}(v)\mbox{,}$$ 
o\`u $v\in O''$, puisque $k_{\phi^{-1}(n+1)}^{-1}(k_{u(i+1)}(w))$ est d\'efini. On a 
$$v = k_{s\lceil (\vert s\vert -1),t}^{-1} k_{u(i+1)}^{} ...k_{u(j-1)}^{-1} 
k_{s\lceil (\vert s\vert -1),t}^{} (z).$$ 
D'o\`u $u(i+1) = s\lceil (\vert s\vert -1),t = u(j-1)$ et donc $i+2<j-1$, ce qui est absurde (on utilise le fait que $z$ et $v$ sont dans $O''$).\bigskip

 L'entier $i$ est donc unique et on peut trouver un ouvert-ferm\'e non vide de 
$D_{k_{\phi^{-1}(n+1)}}$ sur lequel $k_{\phi^{-1}(n+1)}$ co\"\i ncide 
avec une composition des fonctions $k_{\phi^{-1}(0)}$, ..., $k_{\phi^{-1}(n)}$ de la 
forme $k_{\phi^{-1}(j)}H_{m+1}$, o\`u $m<p$. Mais 
ceci est contraire \`a la construction de $k_{\phi^{-1}(n+1)}$. Pour v\'erifier la 
condition (7), on remarque que dans une composition du type de la condition (6) sans 
termes cons\'ecutifs identiques des fonctions $k_{\phi^{-1}(0)}$, ..., $k_{\phi^{-1}(n+1)}$, il y a au plus une fois la fonction $k_{\phi^{-1}(n+1)}$ (comme pr\'ec\'edemment). Une telle composition est donc 
n\'ecessairement de la forme $H_{m+1}$, $H_{m+1}k^{-1}_{\phi^{-1}(i)}k^{}_{\phi^{-1}(n+1)}$, 
$$H_{m+1}k^{-1}_{\phi^{-1}(n+1)}k^{}_{\phi^{-1}(i)}\mbox{,}$$ 
$H_{m+1}k^{-1}_{\phi^{-1}(i)}k^{}_{
\phi^{-1}(n+1)}H_{m'+1}$, $H_{m+1}k^{-1}_{\phi^{-1}(n+1)}k^{}_{\phi^{-1}(i)}H_{m'+1}$, 
$k^{-1}_{\phi^{-1}(i)}k^{}_{\phi^{-1}(n+1)}H_{m+1}$ ou 
$$k^{-1}_{\phi^{-1}(n+1)}k^{}_{\phi^{-1}(i)}H_{m+1}\mbox{,}$$ 
o\`u ${m,~m'\! <\! p}$ et $i\leq n$. Il n'y en a donc qu'un nombre fini.$\hfill\square$

\vfill\eject

 Il reste \`a pouvoir assurer la r\'eduction de la situation d'arriv\'ee au bor\'elien dont 
nous sommes partis. Le lemme qui suit, coupl\'e avec le th\'eor\`eme 3, va le permettre.
 
\begin{lem} Soient $(Z,T,(l_{r,p})_{(r,p)\in\omega^2})$ une situation g\'en\'erale, 
et $(Z,T,(k_{q,p})_{(q,p)\in\omega^2})$ une 
situation d'arriv\'ee telles que pour tout $x$ de $G(k)\cap G(l)$, 
$k[x]\subseteq l[x]$. Alors il existe un 
ensemble $N$, $G_\delta$ de $(Z\times T)\setminus (\bigcup_{r\in\omega} \mbox{Gr}(l_r))$, 
tel que $(Z,T,(k_{q,p})_{(q,p)\in\omega^2},N)$ soit un syst\`eme r\'educteur.\end{lem}

\noindent\bf D\'emonstration.\rm\ Posons $H := G(k)\cap G(l)$ et 
$$N := \{ (x,y)\in Z\times T~/~x\in H~\mbox{et}~y\in\overline{k[x]}\setminus l[x]\}.$$ 
Alors $N$ est clairement $G_\delta$, et on a $N\cap (\bigcup_{r\in\omega} \mbox{Gr}(l_r))=
\emptyset$, ainsi que les conditions (a) et (b) d'un syst\`eme r\'educteur.\bigskip

\noindent $\bullet$ Soit $U$ un ouvert-ferm\'e non vide de $Z$. On choisit $x\in H
\cap U$. Comme $(Z,T,(k_{q,p})_{(q,p)\in\omega^2})$ est une situation d'arriv\'ee, $k[x]$ est 
d\'enombrable sans point isol\'e, 
donc on trouve $y$ dans $\overline{k[x]}\setminus l[x]$, puisque $\overline{k[x]}$ est 
polonais parfait et que $l[x]$ est d\'enombrable. On a alors que $x\in U\cap \Pi_Z[N]$ puisque $(x,y)\in N$. 
D'o\`u la condition (c) d'un syst\`eme r\'educteur.\bigskip

\noindent $\bullet$ Soient $x\in U\cap \Pi_Z[N\cap O]$, $y\in \overline{k[x]}\cap O_x$. On choisit 
des ouverts-ferm\'es $V$ et $W$ tels 
que ${(x,y)\in V\times W}$ et ${V\times W\subseteq O\cap (U\times T)}$. On peut 
trouver $q$ tel que 
$k_q(x)\in W$ et un ouvert-ferm\'e non vide $V'\subseteq k_q^{-1}(W)\cap V$ ; si 
$z\in V'\cap H$, 
$\overline{k[z]}\setminus l[z]$ \'etant dense dans $\overline{k[z]}$, on peut trouver 
$y(z)$ tel que $(z,y(z))\in N\cap (V\times W)$. Donc $U\cap\Pi_Z[N\cap O]$ contient 
$V'\cap H$, qui est non maigre. On a donc montr\'e que pour tout ouvert-ferm\'e 
non vide $U$ de $Z$, $U\cap \Pi_Z[N\cap O]$ est non maigre. Comme $\Pi_Z[N\cap O]$ est 
analytique, on en d\'eduit que $\Pi_Z[N\cap O]$ est comaigre dans $Z$. D'o\`u la condition 
(d) d'un syst\`eme r\'educteur.\bigskip

\noindent $\bullet$ Soient $U$ et $V$ des ouverts-ferm\'es tels que ${N\cap (U\times V)\not=
\emptyset}$, et ${(x,y)\in N\cap (U\times V)}$. Comme $x\in G(k)$ et $y\in V\cap 
\overline{k[x]}$, $k[x]$ n'a pas de point isol\'e et on peut trouver une infinit\'e de $q$ 
tels que $z := k_q(x)\in V$. Comme avant, on voit que $\overline{k[x]}\setminus l[x]$ 
est dense dans $\overline{k[x]}$, donc $z$ est limite de points $z_n\in V\cap
\overline{k[x]}\setminus l[x]$. Donc $(x,z_n)$ est dans ${N\cap (U\times V)}$, et $(x,z)$ 
est dans ${\mbox{Gr}(k_q)\cap\overline{N}\cap (U\times V)}$. Par 
cons\'equent, l'ensemble ${\{q\in\omega~/~\mbox{Gr}(k_q)\cap \overline{N}\cap (U\times V)\not=
\emptyset\}}$ est infini.\bigskip

 Posons $U_{q,p,U,V} := U\cap k_{q,p}^{-1}(V)$. Il est clair que $\mbox{Gr}(k_{q,p})\cap \overline{N}
\cap (U\times V)\subseteq \mbox{Gr}(k_{q,p}\lceil U_{q,p,U,V})$. R\'eciproquement, si $(x,y)\in 
\mbox{Gr}(k_{q,p}\lceil U_{q,p,U,V})$, il faut voir que $(x,y)\in\overline{N}$. $x$ est limite de 
$$(x_n)\subseteq D_{k_{q,p}}\cap H.$$
Comme avant, $k_{q,p}(x_n)$ est limite de 
$(y^n_m)_m\subseteq \overline{k[x_n]}\setminus l[x_n]$, et on peut supposer que 
$$d(k_{q,p}(x_n),y^n_m)<2^{-n-m}.$$
Alors $(x_n,y^n_n)\in N$ et tend vers $(x,y)\in 
\overline{N}$.$\hfill\square$

\vfill\eject

\noindent\bf {\Large 3 Existence d'exemples et synth\`ese des r\'esultats pr\'ec\'edents.}\rm\bigskip 

\noindent\bf Notations.~\rm Soit $(q_n)$ la suite des nombres premiers : 
$q_0 = 2$, $q_1 = 3$, $q_2 = 5$, ... On pose 
$$J : \left\{\!\!
\begin{array}{ll}
\omega^{<\omega}\!\!\!\! 
& \rightarrow\omega \cr 
s 
& \mapsto \left\{\!\!\!\!\!\!
\begin{array}{ll} 
& q_0^{s(0)+1}...q_{\vert s\vert -1}^{s(\vert s\vert -1)+1}~\mbox{si}~s\not=\emptyset\mbox{,}\cr 
& 0~\mbox{sinon.}
\end{array}\right.
\end{array}
\right.$$
$\bullet$ On d\'efinit $A_i := \{1\}\cup \{J(u^\frown 1)~/~u\in\Pi_{p<i}~A_p\}$, puis 
$F := \Pi_{i\in\omega}~A_i$. \bigskip

\noindent $\bullet$ On pose, pour $(s,t)\in\bigcup_{n\in\omega} \omega^n\times\omega^{n+1}$, 
$$D_{f_{s,t}} := \{\alpha\in F~/~\forall~j\leq\vert s\vert ~~\alpha (J[s\lceil j^\frown t
\lceil (j+1)])=1~\mbox{et}~\forall~p<t(j)~~\alpha (J[s\lceil j^\frown t\lceil j^
\frown p])\not=1\}.$$ 
 On d\'efinit ensuite les fonctions 
$$f_{s,t}:\left\{\!\!
\begin{array}{ll}
D_{f_{s,t}}\!\!\!\! 
& \rightarrow F \cr 
\alpha 
& \mapsto\left\{\!\!
\begin{array}{ll} 
\omega\!\!\!\!\!\!\!\!\!\!
& \rightarrow\omega \cr 
q 
& \mapsto \left\{\!\!\!\!\!\!\!\!
\begin{array}{ll} 
& \alpha (q)~\mbox{si}~\forall~j\leq\vert s\vert ~~q\not= J[s\lceil j^\frown t\lceil (j+1)]\mbox{,}\cr\cr 
& J(\alpha\lceil (q\!  +\! 1))~\mbox{si}~\exists~j\leq\vert s\vert ~~q\! =\!  J[s\lceil j^\frown t\lceil (j+1)].
\end{array}
\right.
\end{array}
\right.
\end{array}
\right.$$

\begin{lem} Il existe une situation de d\'epart.\end{lem}

\noindent\bf D\'emonstration.\rm\ Soient $e:\omega\rightarrow\omega^{<\omega}$ 
bijective v\'erifiant $e^{-1}(s) < e^{-1}(s^\frown n)$, et aussi 
$${\Psi : \omega^2 \rightarrow\bigcup_{n\in\omega} \omega^n\times\omega^{n+1}}$$ 
et $\theta : \bigcup_{n\in\omega}\omega^n\times\omega^{n+1} \rightarrow \omega^2$ bijectives 
r\'eciproques l'une de l'autre avec $\Psi_0 (n,p) = e(n)$ et $\theta_0(s,t^\frown m) = e^{-1}(s)$. On pose 
$f_{n,p} := f_{\Psi (n,p)}$. On va voir que $(F,(f_{n,p})_{(n,p)\in\omega^2})$ est une 
situation de d\'epart.\bigskip

\noindent $\bullet$ Il est clair que $A_i$ est une partie finie de 
$\omega$, de cardinal au moins deux. De sorte que $F$ est un compact 
parfait non vide de $\omega^\omega$, donc une copie de $2^\omega$. D'o\`u la condition (a) d'une 
situation g\'en\'erale et la condition (b) d'une situation de d\'epart.\bigskip

\noindent $\bullet$ Les ensembles $D_{f_{s,t}}$ sont des ouverts-ferm\'es de $F$, donc des compacts, 
et $f_{s,t}$ est clairement d\'efinie, injective et continue, donc est un 
hom\'eomorphisme de $D_{f_{s,t}}$ sur son image 
$${\{\alpha\!\in\! F~/~\exists~u\!\in\! \Pi_{i\leq J(s^\frown t)}~A_i~~\forall~i\!
\leq\! J(s^\frown t)~~\alpha (i)\! =\! f_{s,t} (u^\frown 1^\omega )(i)~\mbox{et}~u^\frown 1^\omega\in 
D_{f_{s,t}}\}}.$$
Cette image est ouverte-ferm\'ee dans $F$, d'o\`u la condition (b) d'une situation g\'en\'erale.\bigskip

\noindent $\bullet$ Il est clair que les $D_{f_{s,t^\frown m}}$, pour $m\in \omega$, sont deux \`a deux 
disjoints, donc les $D_{f_{s,t}}$, pour $t\in\omega^{\vert s\vert +1}$, aussi car 
$D_{f_{s,t^\frown m}}\subseteq D_{f_{s\lceil (\vert s\vert -1),t}}$. Leur r\'eunion 
$D_{f_s}$ est dense dans $F$, car si $\tilde t\in \Pi_{i\leq k}~A_i$, on 
peut trouver $t\in \omega^{\vert s\vert +1}$ telle que $\tilde t^\frown 1^\omega \in D_{f_{s,t}}$.

\vfill\eject

 En effet, on peut trouver une suite finie d'entiers $\tilde s$ telle que pour tout 
$j<\vert\tilde s\vert\leq\vert s\vert +1$, on ait $\tilde t(J[s\lceil j^\frown \tilde s\lceil 
(j+1)])=1$ et pour tout 
$p<\tilde s(j)$, $\tilde t(J[s\lceil j^\frown \tilde s\lceil j^\frown p])\not= 1$ ; de plus, on peut exiger 
qu'elle soit de longueur maximale avec ces propri\'et\'es, c'est-\`a-dire que pour tout 
$i\in\omega$, $J[s\lceil \vert\tilde s\vert^\frown \tilde s^\frown i]\geq \vert\tilde t\vert $ 
ou $\tilde t(J[s\lceil \vert\tilde s\vert^\frown \tilde s^\frown i])\not= 1$ ou 
$\vert\tilde s\vert = \vert s\vert +1$. 
Si $\vert\tilde s\vert = \vert s\vert +1$, on peut poser $t := \tilde s$. Sinon, on choisit 
$i$ minimal tel que 
$J[s\lceil \vert\tilde s\vert^\frown \tilde s^\frown i]\geq \vert\tilde t\vert $ et on pose 
$t := \tilde s^\frown i^\frown 0^{\vert s\vert -\vert\tilde s\vert}$. D'o\`u la condition (c) d'une situation 
g\'en\'erale.\bigskip

\noindent $\bullet$ Si $x\in D_{f_{s^\frown n,t^\frown m}}$, $x\in D_{f_{s,t}}$ et on a 
$d(f_{s^\frown n,t^\frown m}(x),f_{s,t}(x))\leq 2^{-J(s^\frown n^\frown t^\frown m)}$. 
Donc si $x$ est dans $D_{f_{s^\frown n}}\cap D_{f_s}$, $\exists (t,m)\in \omega^{\vert s\vert +1}
\times\omega$ tel que $x\in D_{f_{s^\frown n,t^\frown m}}$ et la suite 
$(f_{s^\frown n}(x))_n$ tend vers $f_s(x)$, qui n'est donc pas isol\'e. On a donc la 
condition (d) d'une situation g\'en\'erale, avec $G(f) := \bigcap_{n\in\omega} D_{f_n}$ ; 
d'o\`u la condition (a) d'une situation de d\'epart.\bigskip

\noindent $\bullet$ La condition (c) d'une situation de d\'epart est clairement v\'erifi\'ee. Soit donc 
$(x,y)$ dans 
$$(\bigcap_{n\in\omega} D_{f_n}\times F)\cap 
\overline{\bigcup_{n\in\omega} \mbox{Gr}(f_n)}.$$ 
On peut trouver $x_k$ et $n_k$ tels que $(x,y)$ soit la limite en $k$ de $(x_k,f_{n_k}(x_k))$ : 
$$\forall~p\in\omega~~\exists~k(p)\in\omega~~\forall~k\geq k(p)~~x_k\lceil p = 
x\lceil p~\mbox{et}~f_{n_k}(x_k)\lceil p = y\lceil p.$$
 Montrons que la suite $(f_n)$ est \'equicontinue en $x$. Soit $q \in \omega$. Comme $x$ 
est dans $\bigcap_{n\in\omega} D_{f_n}$, les coordonn\'ees de $x$ modifi\'ees 
sont les $J(s^\frown t^\frown m)$, o\`u $x\in D_{f_{s,t^\frown m}}$, $s,t\in
\omega^{<\omega}$ et $\vert s\vert =\vert t\vert $, $m$ variant dans $\omega$. Par suite, 
les couples $(x\lceil q,f_n(x)\lceil q)$ sont donc en 
nombre fini, \`a $q$ fix\'e. Comme les $f_n$ sont continues, on a l'\'equicontinuit\'e : 
$$\forall~p\in\omega~~\exists~q(p)\in \omega~~z\lceil q(p)=x\lceil q(p)~\Rightarrow~
\forall~n\in\omega~~f_n(z)\lceil p = f_n(x)\lceil p.$$
Posons $y_k := f_{n_k}(x)$. Alors $(x,y_k)\in \bigcup_{n\in\omega} \mbox{Gr}(f_n)$ 
et si $p\in\omega$, soit $k\geq k(\mbox{max}[p,q(p)])$. Alors $x_k\lceil q(p) = 
x\lceil q(p)$, d'o\`u $y\lceil p = f_{n_k}(x_k)\lceil p = f_{n_k}(x)\lceil p = 
y_k\lceil p$. Pour avoir la condition (d) d'une situation de d\'epart, il reste \`a voir que 
$x\not= y$. Ceci r\'esulte du fait que la premi\`ere coordonn\'ee de $x$ \'egale \`a $1$ et 
d'ordre une puissance positive de $2$ est transform\'ee de la m\^eme fa\c con par toutes les 
$f_{n_k}$. \bigskip

\noindent $\bullet$ Si $N_s^2$ rencontre l'un des graphes des $f_n$, il n'y a pas de $1$ 
parmi les coordonn\'ees de $s$ d'ordre une puissance positive de $2$ ; par suite, $N_s^2$ rencontre 
le graphe de $f_0$. Si maintenant $x\not= y$, $x^\frown j~{\cal R}~y^\frown j$ et 
$(N_x\times N_y)\cap \mbox{Gr}(f_s)\not=\emptyset$, il existe $t\prec s$ telle que 
$(N_x\times N_y)\cap \mbox{Gr}(f_t)\not=\emptyset$, avec $f_t$ ne changeant que les coordonn\'ees d'ordre 
inf\'erieur \`a $\vert x\vert $. De plus, toute suite $u$ telle que 
$(N_x\times N_y)\cap \mbox{Gr}(f_u)\not=\emptyset$ v\'erifie $t\prec u$, d'o\`u 
$e^{-1}(t)\leq e^{-1}(u)$ et $\psi (x,y) = e^{-1}(t)$. Comme 
$(N_{x^\frown j}\times N_{y^\frown j})\cap \mbox{Gr}(f_t)\not=\emptyset$, 
$\psi(x^\frown j,y^\frown j)\leq e^{-1}(t)=\psi(x,y)$. D'o\`u la condition (e).(i) d'une 
situation de d\'epart.\bigskip

\noindent $\bullet$ Soient $s$ et $t$ dans $\Pi_{i\leq k}~A_i$ et $u$, $v$ des $\cal T$-cha\^\i nes 
sans r\'ep\'etition de termes telles que $u(0) = v(0) = s$ et $u(\vert u\vert -1) = v(\vert v\vert -1) 
= t$, avec $\vert u\vert \leq\vert v\vert $. On veut montrer que $u=v$ ; on peut supposer que 
$\vert u\vert \geq 2$ et $u(\vert u\vert -2)\not= v(\vert v\vert -2)$.

\vfill\eject

 Posons 
$$w(i) = \left\{\!\!\!\!\!\!\!
\begin{array}{ll}
& u(i)~\mbox{si}~i<\vert u\vert\mbox{,}\cr 
& v(\vert u\vert +\vert v\vert -2-i)~\mbox{si}~\vert u\vert \leq i<\vert u\vert +\vert v\vert -1.
\end{array}\right.$$
Alors $\vert w\vert\geq 3$, $w(0)=w(\vert w\vert-1)=s$, $w(i)\not=w(i+1)$ si $i<\vert w\vert-1$ 
et $w(i)\not=w(i+2)$ si $i<\vert w\vert-2$. Soit $c$ une $\cal T$-cha\^\i ne de longueur minimale 
ayant ces propri\'et\'es, et telle que $l:= \vert c(0)\vert$ soit minimale elle aussi. \bigskip 

 L'argument qui suit a \'et\'e vu dans la preuve du th\'eor\`eme 2.7 de [Le3]. La suite 
$$(c(i)(l-1))_{i<\vert c\vert }$$ 
est non constante, et on trouve $i_1$ minimal tel que 
$c(i_1)(l-1)\not= c(i_1+1)(l-1)$ ; il y a deux cas.\bigskip

 Ou bien $c(i_1)(l-1)<c(i_1+1)(l-1)$, auquel cas comme on a les \'egalit\'es 
$${c(i_1)(l-1)\! =\! c(0)(l-1)\! =\! c(\vert c\vert -1)(l-1)}\mbox{,}$$ 
on trouve $i_2>i_1+1$ minimal tel que l'on ait ${c(i_1+1)(l-1)\not= c(i_2)(l-1)}$. On a  
que $c(i_1)=c(i_2)$, par injectivit\'e de $J$. Donc $i_1=0$ et 
$i_2=\vert c\vert -1$, par minimalit\'e de $\vert c\vert $. Par minimalit\'e encore, $\vert c\vert =3$, 
ce qui constitue la contradiction cherch\'ee (on a $c(i_1+1) = c(i_2-1)$ car il 
existe un unique couple $(s,t)$ tel que $c(i_1)^\frown 1^\omega\in 
D_{f_{s,t}}$, avec $J(s^\frown t) = l-1$ ; par suite, on a la suite d'\'egalit\'es  
$c(i_1+1)^\frown 1^\omega = f_{s,t} (c(i_1)^\frown 1^\omega) = 
f_{s,t} (c(i_2)^\frown 1^\omega) = c(i_2-1)^\frown 1^\omega$).\bigskip

 Ou bien $c(i_1)(l-1)>c(i_1+1)(l-1)$, auquel cas on trouve $i_2>i_1+1$ minimal 
tel que $c(i_2)(l-1) = ... = c(\vert c\vert -1)(l-1)$. On a $c(i_1+1) = c(i_2-1)$, donc 
$c(i_1)=c(i_2)$ comme avant. D'o\`u $i_1=0$ et $i_2=\vert c\vert -1$, par minimalit\'e de 
$\vert c\vert $. Par minimalit\'e encore, $\vert c\vert =3$, ce qui constitue la contradiction 
cherch\'ee. D'o\`u la condition (e).(ii) d'une situation de d\'epart.$\hfill\square$

\begin{thm} Il existe un bor\'elien $B$ de $\omega^\omega\times\omega^\omega$, tel que pour 
tous espaces polonais $X$ et $Y$, et pour tout bor\'elien $A$ de $X\times Y$ dont les coupes 
horizontales et verticales sont d\'enombrables, on a l'\'equivalence entre les 
conditions suivantes :\smallskip

\noindent (a) Le bor\'elien $A$ n'est pas $\mbox{pot}(\bormtwo)$.\smallskip

\noindent (b) Il existe $u : \omega^\omega \rightarrow X$ et 
$v : \omega^\omega \rightarrow Y$, hom\'eomorphismes sur leurs images, tels que l'on ait
$\overline{B}\cap (u\times v)^{-1}(A) = B$.\end{thm}

\noindent\bf D\'emonstration.\rm\ Soit 
$(F,(f_{n,p})_{(n,p)\in\omega^2})$ la situation de d\'epart fournie par 
le lemme 6. L'ensemble ${G(f) = \bigcap_{n\in\omega} D_{f_n}}$ est 
$G_\delta$ dense de $F$, donc polonais parfait de dimension $0$, et $G(f)$ est 
localement non compact car son compl\'ementaire contient l'ensemble dense $D$ des suites 
diff\'erentes de $1$ \`a partir d'un certain rang. On peut donc trouver un hom\'eomorphisme 
$\phi_0 : \omega^\omega \rightarrow G(f)$. On remarque que si $x\in G(f)$ et $n\in\omega$, 
$f_n(x)\notin D$, \`a cause de la condition (c) d'une situation de d\'epart. Soit donc 
$\psi_0 : \omega^\omega\rightarrow F\setminus D$ un 
hom\'eomorphisme. On pose $$B := (\phi_0\times \psi_0)^{-1}(\bigcup_{n\in\omega} Gr[f_n\lceil G(f)]).$$

\vfill\eject

\noindent $\bullet$ Si $A$ est $\mbox{pot}(\bormtwo)$, alors la condition (b) n'est pas 
v\'erifi\'ee, car sinon $B$ serait $\mbox{pot}(\bormtwo)$, donc $\bigcup_{n\in\omega} 
\mbox{Gr}(f_n\lceil G(f))$ aussi. On pourrait donc trouver un $G_\delta$ dense $K$ de $F$ tel que 
pour tout $x$ de $G(f)$, $f[x]\cap K$ soit $G_\delta$ de $K$, donc 
polonais. Mais $\{x\in G(f)~/~x\in \bigcap_{n\in\omega} f_n^{-1}(K)\}$ est $G_\delta$ 
dense de $G(f)$, donc on pourrait trouver $x$ dans $G(f)$ tel que $f[x]$ 
soit polonais, ce qui contredit le fait qu'il soit sans point isol\'e.\bigskip

\noindent $\bullet$ Si $A$ n'est pas $\mbox{pot}(\bormtwo)$, nous allons construire des applications 
$u$ et $v$ v\'erifiant la condition (b). Le lemme 2 fournit un syst\`eme r\'educteur  
$(Z,T,(h_{n,p})_{(n,p)\in\omega^2},M)$ et $u_0\! :\! Z\!\rightarrow \! X$, 
$v_0\! :\! T\! \rightarrow \! Y$ injectives continues tels que 
${\bigcup_{(n,p)\in\omega^2}\mbox{Gr}(h_{n,p})\! \subseteq \! (u_0\! \times \! v_0)^{-1}(A)}$ et 
${M\subseteq (u_0\times v_0)^{-1}(\check A)}$. Par le th\'eor\`eme 3, on trouve une injection 
$\Psi:\omega^2\rightarrow\omega^2$ et des ouverts-ferm\'es 
$D'_{r,p}\subseteq D_{h_{\Psi(r,p)}}$ tels que si ${l_{r,p} := h_{\Psi(r,p)}\lceil D'_{r,p}}$, 
$(Z,T,(l_{r,p})_{(r,p)\in\omega^2})$ soit une situation g\'en\'erale et $\forall~x\in G(l)$, 
$${l[x] \subseteq \{h_{n,p}(x)~/~(n,p)\in\omega^2
~\mbox{et}~x\in D_{h_{n,p}}\}}$$ 
et ${\forall~y\in \overline{l[x]}\setminus l[x]}$, ${(x,y)\in M}$. Par le th\'eor\`eme 4, 
on trouve une situation d'arriv\'ee 
$$(Z,T,(k_{q,p})_{(q,p)\in\omega^2})$$ 
telle que $G(k)\subseteq G(l)$ et pour $x$ dans $G(k)$, on ait $k[x]\subseteq l[x]$. Par le 
lemme 5, on trouve un ensemble $N$, $G_\delta$ de $Z\times T$, tel que 
$(Z,T,(k_{q,p})_{(q,p)\in\omega^2},N)$ soit un syst\`eme r\'educteur et 
$$N\cap (\bigcup_{r\in\omega} \mbox{Gr}(l_r))=\emptyset.$$ 
Par le th\'eor\`eme 3 encore, il existe une injection 
$\Phi :\omega^2\rightarrow \omega^2$ et des ouverts-ferm\'es 
${D_{m,p}\subseteq D_{k_{\Phi(m,p)}}}$ tels que si 
${g_{m,p} := k_{\Phi(m,p)}\lceil D_{m,p}}$, ${(Z,T,(g_{m,p})_{(m,p)\in\omega^2})}$ 
soit une situation d'arriv\'ee et 
pour $x$ dans ${G(g)\subseteq G(k)}$, ${g[x]\subseteq k[x]}$ et 
pour tout $y$ de ${\overline{g[x]}\setminus g[x]}$, 
$(x,y)\in N$. Par le th\'eor\`eme 1, on trouve des injections continues 
${u_1:F\rightarrow G(g)}$ et ${v_1:F\rightarrow T}$ telles que pour $(x,y)$ dans 
${\bigcup_{n\in\omega} \mbox{Gr}(f_n)}$ on ait ${(u_1(x),v_1(y))\! \in\! \bigcup_{m\in\omega} 
\mbox{Gr}(g_m)}$ 
et ${\forall~(x,y)\! \in \! (\bigcap_{n\in\omega} D_{f_n}\! \times\!  T)\! \cap\! 
\overline{\bigcup_{n\in\omega} \mbox{Gr}(f_n)}
\setminus (\bigcup_{n\in\omega} \mbox{Gr}(f_n))}$, on ait l'appartenance de $v_1(y)$ \`a 
${\overline{g[u_1(x)]}\setminus g[u_1(x)]}$.\bigskip

 On pose alors $u:= u_0\lceil G(g)~\circ ~u_1\lceil G(f)~\circ ~\phi_0$ et 
$v := v_0~\circ ~v_1\lceil (F\setminus D)~\circ ~\psi_0$. Comme $u_0\lceil G(g)~\circ ~u_1$ et 
$v_0~\circ ~v_1$ sont 
des hom\'eomorphismes sur leurs images, $u$ et $v$ aussi. Si ${(x,y)\in B}$, 
${(\phi_0(x),\psi_0(y))\in \bigcup_{n\in\omega} \mbox{Gr}(f_n\lceil G(f))}$ donc 
${(u_1[\phi_0(x)],v_1[\psi_0(y)])\in \bigcup_{m\in\omega} \mbox{Gr}(g_m)}$ et aussi  
$${u_1[\phi_0(x)]\in G(g)\subseteq G(k) \subseteq G(l)}.$$
Par suite, on a 
$$\begin{array}{ll}
v_1[\psi_0(y)] 
& \in g[u_1[\phi_0(x)]]\cr 
& \subseteq k[u_1[\phi_0(x)]]\cr 
& \subseteq l[u_1[\phi_0(x)]]\cr 
& \subseteq \{h_{n,p}(u_1[\phi_0(x)])~/~(n,p)\in \omega^2~\mbox{et}~u_1[\phi_0(x)]\in D_{h_{n,p}}\}.
\end{array}$$
Donc $(u(x),v(y))\in A$.

\vfill\eject

 Si $(x,y)\in\overline{B}\setminus B$, 
$(\phi_0(x),\psi_0(y))\in \overline{\bigcup_{n\in\omega} \mbox{Gr}(f_n\lceil G(f))}\setminus 
(\bigcup_{n\in\omega} Gr[f_n\lceil G(f)])$, et comme $\phi_0(x)\in G(f)
\subseteq \bigcap_{n\in\omega} D_{f_n}$, $(\phi_0(x),\psi_0(y))\! \in \! \overline{
\bigcup_{n\in\omega} \mbox{Gr}(f_n)}\setminus (\bigcup_{n\in\omega} \mbox{Gr}(f_n))$. Par suite, 
$v_1[\psi_0(y)]$ est dans $\overline{g[u_1[\phi_0(x)]]}\setminus 
g[u_1[\phi_0(x)]]$. Comme ${u_1[\phi_0(x)]}$ est dans $G(g)$, on a que 
$(u_1[\phi_0(x)],v_1[\psi_0(y)])$ est dans $N$, et donc que $v_1[\psi_0(y)]$ appartient \`a 
$\overline{l[u_1[\phi_0(x)]]}\setminus l[u_1[\phi_0(x)]]$. Donc 
$(u_1[\phi_0(x)],v_1[\psi_0(y)])$ est dans $M$ et $(u(x),v(y))\notin A$.$\hfill\square$\bigskip

 En analysant cette d\'emonstration, on obtient d'autres caract\'erisations des bor\'eliens 
\`a coupes d\'enombrables n'\'etant pas $\mbox{pot}(\bormtwo)$. Le corollaire qui suit est \`a 
rapprocher du th\'eor\`eme 2.11 de [Le2].

\begin{cor} Soient $X$ et $Y$ des espaces polonais, et $A$ un bor\'elien de 
$X\times Y$ dont les coupes horizontales et verticales sont d\'enombrables. Les 
conditions suivantes sont \'equivalentes :\smallskip

\noindent (a) Le bor\'elien $A$ n'est pas $\mbox{pot}(\bormtwo)$.\smallskip

\noindent (b) Il existe une situation g\'en\'erale $(Z,T,(g_{m,p})_{(m,p)\in\omega^2})$ et 
$i:Z\rightarrow X$, $j:T\rightarrow Y$ injectives continues telles que pour 
tout $x$ de $G(g)$, on ait $\overline{g[x]}\cap (i\times j)^{-1}(A)_x = g[x]$.\smallskip

\noindent (c) Il existe une situation d'arriv\'ee $(Z,T,(g_{m,p})_{(m,p)\in\omega^2})$ et 
$i:Z\rightarrow X$, $j:T\rightarrow Y$ injectives continues telles que pour 
tout $x$ de $G(g)$, on ait $\overline{g[x]}\cap (i\times j)^{-1}(A)_x = g[x]$.\smallskip

\noindent (d) Il existe une situation de d\'epart $(F,(f_{n,p})_{(n,p)\in\omega^2})$ et ${u:F\!\rightarrow\! X}$, ${v:F\rightarrow Y}$ injectives continues telles que 
$${\overline{\bigcup_{n\in\omega} \mbox{Gr}(f_n)}\cap 
(G(f)\times F)\cap (u\times v)^{-1}(A) = \bigcup_{n\in\omega} \mbox{Gr}(f_n\lceil G(f))}.$$\end{cor}

\noindent\bf D\'emonstration.\rm\ Il suffit de relire la preuve du th\'eor\`eme 7. Pour 
l'\'equivalence de (a) et (d), on prend $u := u_0\lceil G(g)~\circ ~u_1$ et 
$v := v_0~\circ ~v_1$. Pour l'\'equivalence de (a) avec (b) et (c), on prend $i := u_0$ 
et $j := v_0$, de sorte que $i$ et $j$ correspondent simplement \`a un changement de 
topologie. Ces \'equivalences viennent du fait que ${\bigcup_{m\in\omega} \mbox{Gr}(g_m\lceil 
G(g)) = \{(x,y)\in G(g)\times T~/~y\in\overline{g[x]}\}\cap (i\times j)^{-1}(A)}$.$\hfill\square$

\begin{cor} Soit $\Gamma$ une classe de Wadge non stable par passage au 
compl\'ementaire. Alors il existe un bor\'elien $B_\Gamma$ de $\omega^\omega\times\omega^\omega$ et un ferm\'e $F_\Gamma$ contenant $B_\Gamma$, tels que pour tous espaces polonais $X$ et $Y$, et pour tout bor\'elien $A$ de $X\times Y$ ayant ses coupes horizontales et verticales d\'enombrables, on a  
l'\'equivalence entre les conditions suivantes :\smallskip

\noindent (a) Le bor\'elien $A$ n'est pas $\mbox{pot}(\Gamma )$.\smallskip

\noindent (b) Il existe des fonctions continues $u : \omega^\omega\rightarrow X$ et 
$v : \omega^\omega\rightarrow Y$ telles que $F_\Gamma\cap (u\times v)^{-1}(A) = B_\Gamma$.\end{cor}

\noindent\bf D\'emonstration.\rm\ Si $\Gamma = \bormtwo$, on applique le th\'eor\`eme 7. Si 
$\Gamma = D_\xi (\boraone)$ ou $\check D_\xi (\boraone)$, on applique les th\'eor\`emes 3.5 et 
3.6 de [Le3] et on utilise l'existence d'une r\'etraction continue de 
$\omega^\omega$ sur $2^\omega$. Sinon, $\Gamma$ contient $\boratwo$, donc $A$ est 
$\mbox{pot}(\Gamma)$. Il suffit alors de prendre $B_\Gamma := (\omega^\omega\times\omega^\omega)
\setminus B_{\bormtwo}$ et $F_\Gamma := \overline{B_\Gamma}$, puisque $B_{\bormtwo}
\notin \mbox{pot}(\bormtwo)$, par le th\'eor\`eme 7.$\hfill\square$\bigskip

\noindent \bf Remarque.\rm ~$B_{\bormtwo}$ \'etant r\'eunion d\'enombrable de graphes de fonctions continues 
est ${\boratwo\setminus \mbox{pot}(\bormtwo)}$. Si $\Gamma\subseteq \bortwo$ et $\check \Gamma$ est 
stable par intersection avec les ferm\'es (c'est-\`a-dire si $\Gamma = D_\xi (\boraone)$ 
avec $\xi$ impair ou $\Gamma = \check D_\xi (\boraone)$ 
avec $\xi$ pair), on a aussi $B_\Gamma \in \check \Gamma\setminus \mbox{pot}(\Gamma)$. En effet, 
on applique le th\'eor\`eme B \`a $A\in \check \Gamma\setminus \mbox{pot}(\Gamma)$ (qui existe par 
le th\'eor\`eme 3.3 de [Le1]) pour voir que $B_\Gamma = \overline{A_\xi}\setminus 
A_\xi \in\check \Gamma$. On en d\'eduit que $B_{\check\Gamma} = A_\xi \in D_{\xi+1} (\boraone)$ 
si $\xi$ est impair et que $B_{\check\Gamma} \in\check D_{\xi+1} (\boraone)$ 
si $\xi$ est pair. On n'a pas mieux en g\'en\'eral (cf [Le3] pour $\xi = 1$ : 
$B_{\bormone} \in D_2 (\boraone) \setminus \mbox{pot}(\check D_2 (\boraone))$).

\vfill\eject

\noindent\bf {\Large 4 Une limite du r\'esultat principal.}\bigskip\rm

 On peut observer un ph\'enom\`ene analogue \`a celui d\'ecrit dans la section 2.C de 
[Le3], c'est-\`a-dire que dans le th\'eor\`eme 7, en supposant seulement $A$ \`a 
coupes verticales d\'enombrables, on a une incompatibilit\'e avec l'existence des 
injections $u$ et $v$.\bigskip

\noindent\bf D\'efinition.\it\ Une suite $(h_s)_{s\in (\omega\setminus\{ 0\} )^{<\omega}}$ de fonctions partielles de $2^\omega$ dans $2^\omega$ est une $bonne\ suite$ si\smallskip

\noindent (a) Le domaine $D_{h_s}$ de $h_s$ est un ouvert dense de $2^\omega$.\smallskip

\noindent (b) Les fonctions $h_s$ sont continues et ouvertes.\smallskip

\noindent (c) Si $x\in \bigcap_{s\in (\omega\setminus\{ 0\} )^{<\omega}} D_{h_s}$ et 
$s\in (\omega\setminus\{ 0\} )^{<\omega}$, $\displaystyle \lim_{k\rightarrow \infty} 
{h_{s^\frown k}(x)} = h_s(x)$.\smallskip

\noindent (d) Si $s,~t\in (\omega\setminus\{ 0\} )^{<\omega}$ et $s\not= t$, 
$\{x\in D_{h_s}\cap D_{h_t}~/~h_s(x)\not= h_t(x)\}$ est dense dans $2^\omega$.\rm\bigskip

\bf\noindent Exemple.\rm~Soit $J:\omega^{<\omega}\rightarrow\omega$ l'injection 
d\'efinie dans la section 3. On pose 
$$h_s\! :\!  \left\{\!\! 
\begin{array}{ll} 
2^\omega\!\!\!\! 
& \! \rightarrow\!  2^\omega \cr 
x 
& \! \mapsto\!  \left\{\!\! 
\begin{array}{ll} 
\omega\!\!\!\! 
& \! \rightarrow\!  2\cr 
k 
& \! \mapsto\!  \left\{\!\!\!\!\!\!\! 
\begin{array}{ll} 
& x(J(s\lceil i)q-{J(s\lceil i)\over q_{i-1}}-1)~\mbox{si}~\left\{\!\!\!\!\!\! 
\begin{array}{ll} 
& k={J(s\lceil i)\over q_{i-1}}q-1\cr 
& \mbox{et}\cr 
& i\! =\! \mbox{max}\{1\! \leq\!  j\! \leq\! \vert s\vert\!\! ~/~\! k\! \equiv \! -1~({J(s\lceil j)\over q_{j-1}})\}\mbox{,}
\end{array}\right.\cr\cr 
& x(k)~\mbox{si}~\forall~1\leq i\leq\vert s\vert ~~k\not\equiv -1~({J(s\lceil i)\over q_{i-1}}).
\end{array}
\right.
\end{array}
\right.
\end{array}\right.$$
Alors $(h_s)_{s\in (\omega\setminus\{ 0\} )^{<\omega}}$ est une bonne suite. En effet, 
les conditions (a) et (b) sont clairement r\'ealis\'ees. Le plus petit entier $n$ tel que 
$h_{s^\frown k}(x)(n) \not= h_s(x)(n)$, s'il existe, est sup\'erieur ou \'egal \`a 
${J(s^\frown k)\over q_{\vert s\vert }}-1$, qui tend vers l'infini avec $k$. D'o\`u la 
condition (c). Soit $u\in 2^{<\omega}$ ; on cherche $x\in2^\omega$ tel que 
$h_s(u^\frown x)\not= h_t(u^\frown x)$.\bigskip

 Ou bien on trouve $m<\mbox{min}(\vert s\vert ,\vert t\vert )$ tel que $s(m)\not= t(m)$ et $s\lceil m 
=t\lceil m$, avec par exemple $s(m)<t(m)$. Posons alors 
$k_n := {J(t\lceil m+1)\over q_{m}}q_{m+2}^n-1$. On a 
$$\begin{array}{ll} 
h_s(u^\frown x)(k_n)\!\!\!\! 
& = u^\frown x(J(s\lceil m+1)q_m^{t(m)-s(m)-1}q_{m+2}^n-{J(s\lceil m+1)\over q_{m}}-1)\mbox{,}\cr 
h_t(u^\frown x)(k_n)\!\!\!\! 
& = u^\frown x(J(t\lceil m+1)q_{m+2}^n-{J(t\lceil m+1)\over q_{m}}-1).
\end{array}$$
On choisit $n\in\omega$ tel que les num\'eros des coordonn\'ees ci-dessus soient 
sup\'erieurs \`a $\vert u\vert $. Une simplification par $q_{m}^{s(m)}$ montre que ces deux 
num\'eros sont diff\'erents. D'o\`u l'existence de $x$.\bigskip

 Ou bien par exemple $s$ est un d\'ebut strict de $t$. Posons 
$k_n := {J(t)\over q_{\vert t\vert -1}}q_{\vert t\vert +1}^n-1$. On a 
$$\begin{array}{ll} 
h_s(u^\frown x)(k_n)\!\!\!\! 
& = u^\frown x(J(s)q_{\vert s\vert -1}q_{\vert s\vert }^{t(\vert s\vert )+1}...q_{\vert t\vert -2}^{t(\vert t\vert -2)+1}q_{\vert t\vert -1}^{t(\vert t\vert -1)}q_{\vert t\vert +1}^n-{J(s)\over q_{\vert s\vert -1}}-1)\mbox{,}\cr
h_t(u^\frown x)(k_n)\!\!\!\! 
& = u^\frown x(J(t)q_{\vert t\vert +1}^n-{J(t)\over q_{\vert t\vert -1}}-1).
\end{array}$$
On conclut comme avant, avec simplification par $q_{\vert s\vert -1}^{s(\vert s\vert -1)}$.

\vfill\eject

\begin{lem} Soit $(h_s)_{s\in (\omega\setminus\{ 0\} )^{<\omega}}$ une bonne 
suite, et $G$ un $G_\delta$ dense de $2^\omega$ inclus dans 
$$\bigcap_{s\in (\omega\setminus\{ 0\} )^{<\omega}} D_{h_s}\cap\bigcap_{s\not= t} 
\{x\in D_{h_s}\cap D_{h_t}~/~h_s(x)\not= h_t(x)\}.$$ 
Alors $\bigcup_{s\in (\omega\setminus\{ 0\} )^{<\omega}} \mbox{Gr}(h_s\lceil G)$ n'est pas 
$\mbox{pot}(G_\delta)$.\end{lem}

\noindent\bf D\'emonstration.\rm\ Elle est identique \`a celle du deuxi\`eme point de la preuve du 
th\'eor\`eme 7.$\hfill\square$

\begin{lem} Soit $(h_s)_{s\in (\omega\setminus\{ 0\} )^{<\omega}}$ une bonne 
suite. Alors il existe une bonne suite 
$(l_s)_{s\in (\omega\setminus\{ 0\} )^{<\omega}}$, une suite 
$(V_{s,t})_{(s,t)\in \bigcup_{n\in\omega} (\omega\setminus\{ 0\})^n\times 
(\omega\setminus\{ 0\})^{n+1}}$ d'ouverts non vides de $2^\omega$ et 
$$\phi : \bigcup_{n\in\omega} (\omega\setminus\{ 0\})^n\times (\omega\setminus\{ 0\})^
{n+1} \rightarrow (\omega\setminus\{ 0\} )^{<\omega}$$ 
telles que\smallskip

\noindent (a) $D_{l_s} = \bigcup_{t\in\omega^{\vert s\vert +1},~\mbox{disj.}}~ V_{s,t}$ et 
$V_{s^\frown m,t^\frown n}\subseteq V_{s,t}$.\smallskip

\noindent (b) Pour tout $x$ de $V_{s,t}$ on a $l_s(x) = h_{\phi (s,t)}(x)$.\smallskip

\noindent (c) La suite $\phi (s,t)$ est un d\'ebut strict de $\phi (s^\frown m,t^\frown n)$.\smallskip

\noindent (d) On a $D_{l_{s^\frown k}}\subseteq D_{l_s}\cap \bigcap_{j<k} 
D_{l_{s^\frown j}}$ et pour tout $x$ dans $D_{l_{s^\frown k}}$ et tout $i\leq \vert s\vert $\mbox{, }
$l_{s^\frown k}(x) \not= l_{s\lceil i}(x)$.\smallskip

\noindent (e) Pour tout $x$ de $D_{l_{s^\frown k}}$ on a 
$d(l_{s^\frown k}(x),l_s(x)) <\varepsilon (s^\frown k,x)$, o\`u $\varepsilon (s^\frown k,x)$ est par d\'efinition 
$\mbox{min}[2^{-k}\mbox{, min}_{i<\vert s\vert }~{1\over 4} d(l_{s\lceil i+1}(x)\mbox{,}l_{s\lceil i}(x))
\mbox{, min}_{j<k}~{1\over 4} d(l_{s^\frown j}(x)\mbox{,}l_{s}(x))]$.\end{lem}

\noindent\bf D\'emonstration.\rm\ On commence par poser 
$${V_{\emptyset ,n} := \{x\in D_{h_\emptyset}~/~x(n)=1~\mbox{et}~\forall~p<n~~x(p)=0
\}},~~\phi (\emptyset ,n) := \emptyset .$$
Admettons avoir construit $l_y$ pour $y$ d\'ebutant $s$, et pour $y = s^\frown j$ avec 
$j<k$. On va construire $l_{s^\frown k}$. Soit $t\in \omega^{\vert s\vert +1}$ et $(x_n)$ 
une suite dense de l'ensemble suivant : 
$$V_{s,t}\cap D_{l_s}\cap \bigcap_{j<k} D_{l_{s^\frown j}}\cap 
\bigcap_{w\in (\omega\setminus\{ 0\} )^{<\omega}} D_{h_w}\cap\bigcap_{u\not= v} 
\{x\in D_{h_u}\cap D_{h_v}~/~h_u(x)\not= h_v(x)\}.$$ 
Par hypoth\`ese de r\'ecurrence, la fonction 
$\varepsilon(s^\frown k,.) : D_{l_s}\cap \bigcap_{j<k} D_{l_{s^\frown j}}\rightarrow\mathbb{R}_+^*$ est 
d\'efinie et continue. On peut donc trouver un voisinage ouvert $V$ de $x_0$ tel que pour tout $x$ de $V$, on ait 
$$\varepsilon(s^\frown k,x) > {1\over 2} \varepsilon(s^\frown k,x_0).$$ 
Par ailleurs, la condition (c) d'une bonne suite fournit $r_0\geq k$ tel que 
$$d(h_{\phi (s,t)^\frown r_0}(x_0),h_{\phi (s,t)}(x_0)) < {1\over 2} \varepsilon(s^\frown k,x_0).$$ 
Par continuit\'e de $h_{\phi (s,t)^\frown r_0}$ et $h_{\phi (s,t)}$, on trouve un voisinage ouvert $W$ de 
$x_0$ tel que pour tout $x$ de $W$, on ait 
$d(h_{\phi (s,t)^\frown r_0}(x),h_{\phi (s,t)}(x)) < {1\over 2} \varepsilon(s^\frown k,x_0)$ et $x$ soit dans l'ensemble suivant : 
$$V\cap V_{s,t}\cap D_{l_s}\cap 
\bigcap_{j<k} D_{l_{s^\frown j}}\cap \{z\in D_{h_{\phi (s,t)^\frown r_0}}~/~\forall~i
\leq\vert s\vert ~~h_{\phi (s,t)^\frown r_0}(z)\not= l_{s\lceil i}(z)\}.$$ 

\vfill\eject

 On pose $\phi (s^\frown k, t^\frown 0) := \phi (s,t)^\frown r_0$. On choisit 
$V_{s^\frown k, t^\frown 0}$ contenant $x_0$ dans $W$ tel que 
$$V_{s,t}\setminus\overline{V_{s^\frown k, t^\frown 0}} \not=\emptyset .$$ 
Puis on recommence ceci en rempla\c cant $x_0$ par $x_n$, avec $n$ minimal tel que $x_n\notin 
\overline{V_{s^\frown k, t^\frown 0}}$. Le choix de $V_{s^\frown k, t^\frown 1}$ se 
fait dans $V_{s,t}\setminus\overline{V_{s^\frown k, t^\frown 0}}$, avec $x_n\in 
V_{s^\frown k, t^\frown 1}$ et $V_{s,t}\setminus\overline{V_{s^\frown k, t^\frown 0}
\cup V_{s^\frown k, t^\frown 1}} \not=\emptyset$. En it\'erant cette construction, on 
construit $l_{s^\frown k}$, $V_{s^\frown k, t^\frown m}$ et $\phi (s^\frown k, 
t^\frown m)$ v\'erifiant les propri\'et\'es demand\'ees. En effet, on a 
$V_{s,t}\subseteq \overline{\bigcup_{m\in\omega} V_{s^\frown k, t^\frown m}}$, donc 
$$\bigcup_{t\in\omega^{\vert s\vert +1}}~ V_{s,t} \subseteq 
\overline{\bigcup_{t\in\omega^{\vert s\vert +1},~m\in\omega} V_{s^\frown k, t^\frown m}} = 
\overline{D_{l_{s^\frown k}}}.$$ 
La seule chose restant \`a v\'erifier est la condition (d) d'une bonne suite. Elle r\'esultera imm\'ediatement du lemme qui suit.$\hfill\square$

\begin{lem} Soit $(l_s)_{s\in (\omega\setminus\{ 0\} )^{<\omega}}$ une suite 
de fonctions partielles de $2^\omega$ dans $2^\omega$ v\'erifiant les conditions (d) et (e) 
du lemme 11. Soient $s,t\in (\omega\setminus\{ 0\} )^{<\omega}$ et $i<\mbox{min}~
(\vert s\vert ,~\vert t\vert )$ tel que ${s(i) < t(i)}$ et $s\lceil i = t\lceil i$. Alors pour tout 
$x$ de $D_{l_s}\cap D_{l_t}$ on a 
$$d(l_s(x), l_t(x)) \geq {1\over 3} d(l_{s\lceil i+1}(x), l_{s\lceil i}(x)).$$\end{lem}

\noindent\bf D\'emonstration.\rm\ On a 
$$\begin{array}{ll}
d(l_{s\lceil i+1}(x), l_{s}(x))\!\!\!\! & \leq d(l_{s\lceil i+1}(x),l_{s\lceil i+2}(x))+...+
d(l_{s\lceil \vert s\vert -1}(x), l_{s}(x)) \cr 
& \leq d(l_{s\lceil i+1}(x), l_{s\lceil i+2}(x)) (1+{1\over 4}+...+{1\over {4^{\vert s\vert -i-2}}}) \cr 
& \leq {4\over 3} d(l_{s\lceil i+1}(x), l_{s\lceil i+2}(x))~(\mbox{si}~\vert s\vert \geq i+2) \cr 
& \leq {1\over 3} d(l_{s\lceil i+1}(x), l_{s\lceil i}(x))~(\mbox{m\^eme~si}~\vert s\vert = i+1).
\end{array}$$
 De m\^eme, $d(l_{t\lceil i+1}(x), l_{t}(x)) \leq {1\over 3} 
d(l_{t\lceil i+1}(x), l_{s\lceil i}(x))$, car $s\lceil i = t\lceil i$. Par ailleurs on a 
$$\begin{array}{ll}
d(l_{t\lceil i+1}(x), l_{s\lceil i}(x))\!\!\!\! 
& \leq {1\over 4} d(l_{s\lceil i^\frown t(i)-1}(x), l_{s\lceil i}(x)) \leq ...\cr 
& \leq {1\over {4^{t(i)-s(i)}}} d(l_{s\lceil i+1}(x), l_{s\lceil i}(x)) \cr 
& \leq {1\over 4} d(l_{s\lceil i+1}(x), l_{s\lceil i}(x)).
\end{array}$$ 
D'o\`u 
$$\begin{array}{ll}
d(l_{s\lceil i+1}(x), l_{t\lceil i+1}(x))\!\!\!\! 
& \geq d(l_{s\lceil i+1}(x), l_{s\lceil i}(x)) - d(l_{t\lceil i+1}(x), l_{s\lceil i}(x))\cr 
& \geq {3\over 4} d(l_{s\lceil i+1}(x), l_{s\lceil i}(x)).
\end{array}$$
On a aussi $d(l_{t\lceil i+1}(x), l_{t}(x))\leq {1\over 12} 
d(l_{s\lceil i+1}(x), l_{s\lceil i}(x))$. 
D'o\`u 
$$\begin{array}{ll} 
d(l_s(x), l_t(x))\!\!\!\! 
& \geq d(l_{s\lceil i+1}(x), l_{t\lceil i+1}(x))-d(l_{s\lceil i+1}(x), l_{s}(x))-d(l_{t\lceil i+1}(x), l_{t}(x)) \cr 
& \geq d(l_{s\lceil i+1}(x), l_{s\lceil i}(x)) ({3\over 4}-{1\over 3}-{1\over 12}) \cr 
& \geq {1\over 3} d(l_{s\lceil i+1}(x), l_{s\lceil i}(x)).
\end{array}$$
Ceci termine la preuve.$\hfill\square$\bigskip

\noindent\bf Remarque.\rm ~On a en fait montr\'e que pour tout $x$ de $D_{l_s}\cap D_{l_t}$, 
$l_s(x)\not= l_t(x)$ si $s\not= t$.

\vfill\eject

\begin{lem} Soit $(l_s)_{s\in (\omega\setminus\{ 0\} )^{<\omega}}$ une suite 
de fonctions partielles de $2^\omega$ dans $2^\omega$ v\'erifiant les conditions du lemme 
11. Soient $s\in (\omega\setminus\{ 0\} )^{<\omega}\setminus \{\emptyset\}$, $u\in 2^{<\omega}$ et 
$\alpha\in N_u\cap\bigcap_{w\in (\omega\setminus\{ 0\} )^{<\omega}} D_{l_w}$. Alors 
il existe $\varepsilon (s,u,\alpha )>0$ et $v(s,u,\alpha )\in 2^{<\omega}$ tels que 
$u\prec v(s,u,\alpha )\prec \alpha$ et pour tout $t\in (\omega\setminus\{ 0\} )^
{<\omega}$, pour tout $x\in N_{v(s,u,\alpha )}\cap\bigcap_{w\in (\omega\setminus
\{ 0\} )^{<\omega}} D_{l_w}$, on ait l'implication 
$$d(l_s(x), l_t(x)) < \varepsilon (s,u,\alpha )~\Rightarrow~s~\mbox{et}~t~\mbox{sont~
compatibles.}$$\end{lem}

\noindent\bf D\'emonstration.\rm\ Posons, pour $y\in (\omega\setminus\{0\})^{\vert s\vert }$, 
$\eta (y,x) := d(l_{y}(x),l_{s\lceil \vert s\vert -1}(x))$. Alors la fonction $\eta (y,.)
 : D_{l_{y}}\cap D_{l_{s\lceil \vert s\vert -1}} \rightarrow\mathbb{R}_+^*$ 
est d\'efinie et continue, donc on peut trouver $v(s,u,\alpha )\in 2^{<\omega}$ telle 
que $u\prec v(s,u,\alpha )\prec \alpha$ et pour tout $x\in N_{v(s,u,\alpha )}\cap
\bigcap_{w\in (\omega\setminus\{ 0\} )^{<\omega}} D_{l_w}$, on ait $\eta (y,x)> 
{3\over 4} \eta (y,\alpha )$ et $\varepsilon (y,x)> {3\over 4} \varepsilon (y,\alpha )$, ceci pour 
$y\in (\omega\setminus\{0\})^{\vert s\vert }$ v\'erifiant 
$$y=s~\mbox{ou}~(\exists !~i<\vert s\vert ~~y(i)\not= s(i)~\mbox{et}~\exists~i<\vert s\vert ~~y(i)<s(i)).$$
Notons $Y$ l'ensemble de ces suites $y$. On pose $\varepsilon (s,u,\alpha ) := 
\mbox{min}_{y\in Y}~\mbox{min}({1\over 4} \eta (y,\alpha ), \varepsilon (y,\alpha ))$. Soient 
$t\in (\omega\setminus\{ 0\} )^
{<\omega}$ et $x\in N_{v(s,u,\alpha )}\cap\bigcap_{w\in (\omega\setminus
\{ 0\} )^{<\omega}} D_{l_w}$ tels que $d(l_s(x), l_t(x)) < \varepsilon (s,u,\alpha )$. On 
raisonne par l'absurde, ce qui fournit $i<\mbox{min}(\vert s\vert ,\vert t\vert )$ minimal tel que 
$s(i)\not= t(i)$.\bigskip

 Ou bien $s(i)\! <\! t(i)$ ; par le lemme 12, on a $d(l_s(x), l_t(x)) \geq {1\over 3} 
d(l_{s\lceil i+1}(x), l_{s\lceil i}(x))$. Si $\vert s\vert \geq i+2$, on a donc 
$d(l_s(x), l_t(x)) \geq {4\over 3}\varepsilon (s,x) >\varepsilon (s,\alpha ) \geq \varepsilon (s,u,\alpha )$. 
Si $\vert s\vert  = i+1$, on a 
$$d(l_s(x), l_t(x)) \geq {1\over 3}\eta (s,x) >{1\over 4}\eta (s,\alpha ) \geq \varepsilon (s,u,\alpha ).$$ 
Dans tous les cas, on a une contradiction.\bigskip

 Ou bien $t(i)\! <\! s(i)$ ; par le lemme 12, on a $d(l_s(x), l_t(x))\!\geq\! {1\over 3} 
d(l_{t\lceil i+1}(x), l_{t\lceil i}(x))$. Si $\vert s\vert\!\geq\! i\! +\! 2$, 
$d(l_s(x), l_t(x))\!\geq\! {4\over 3}\varepsilon (s\lceil i^\frown t(i)^\frown s(i\! +\! 1),x)\! >\!
\varepsilon (s\lceil i^\frown t(i)^\frown s(i\! +\! 1),\alpha )\!\geq\!\varepsilon (s,u,\alpha )$. 
Si ${\vert s\vert\!  =\! i\! +\! 1}$, on a $d(l_s(x), l_t(x)) \geq {1\over 3}\eta 
(s\lceil i^\frown t(i),x) >{1\over 4}\eta 
(s\lceil i^\frown t(i),\alpha ) \geq \varepsilon (s,u,\alpha )$. L\`a encore, on a une 
contradiction dans les deux cas.$\hfill\square$

\begin{lem} Soit $(l_s)_{s\in (\omega\setminus\{ 0\} )^{<\omega}}$ la bonne 
suite fournie par le lemme 11, associ\'ee \`a $(h_s)_{s\in (\omega\setminus\{ 0\} )^{<\omega}}$ de 
l'exemple. Supposons que $s_i\in (\omega\setminus\{ 0\})^{<
\omega}$, pour $i\in 4$, et que $t_i\in (\omega\setminus\{ 0\})^{\vert s_i\vert+1}$ 
v\'erifient les conditions suivantes :\smallskip

\noindent (a) $\emptyset\not= s_i {\displaystyle \prec_{\not=}} {s_{i+1}}$.\smallskip

\noindent (b) $\emptyset\not= D_{s_i}\subseteq V_{s_i,t_i}$.\smallskip

\noindent (c) $D_{s_{i+1}}\subseteq D_{s_i}$.\smallskip

\noindent (d) $\mbox{Im}(l'_{i+1})\subseteq \mbox{Im}(l'_{i})$, o\`u $l'_i := l_{s_i}\lceil D_{s_i}$.\smallskip
 
Alors $l'_0$ ou $l'_1$ n'est pas injective.\end{lem}

\noindent\bf D\'emonstration.\rm\ Raisonnons par l'absurde. Comme $D_{s_i}\subseteq V_{s_i,t_i}$, 
on a, pour tout $x$ de $D_{s_i}$, 
$$l'_i(x) = l_{s_i}(x) = h_{\phi (s_i,t_i)}(x).$$ 

\vfill\eject

 Comme $D_{s_{i+1}}\subseteq D_{s_i}$, on a $D_{s_{i+1}}\subseteq V_{s_{i+1},t_{i+1}}\cap 
V_{s_i,t_i}$, donc cette intersection est non vide. Comme $V_{s^\frown m,t^\frown n} 
\subseteq V_{s,t}$ et $s_i {\displaystyle \prec_{\not=}} {s_{i+1}}$, $t_i 
{\displaystyle \prec_{\not=}} {t_{i+1}}$, par disjonction de 
$(V_{s,t})_{t\in \omega^{\vert s\vert +1}}$. On a donc $\phi (s_i,t_i) 
{\displaystyle \prec_{\not=}} \phi (s_{i+1},t_{i+1})$. Posons donc 
$$u_i := \phi (s_i,t_i),~~D_{u_i} := D_{s_i},~~h'_i := h_{u_i}\lceil D_{u_i}.$$
On a $\emptyset\not= u_i {\displaystyle \prec_{\not=}} {u_{i+1}}$, $\emptyset\not= 
D_{u_i}$ et $D_{u_{i+1}}\subseteq D_{u_i}$. Les fonctions $h'_i$ et $l'_i$ sont \'egales, donc 
$$\mbox{Im}(h'_{i+1})\subseteq \mbox{Im}(h'_{i})$$ 
et $h'_0$, $h'_1$ sont injectives. De plus, on a, puisque $s_2\not= s_3$, 
$h'_2(x) = l_{s_2}(x) \not= l_{s_3}(x) = h'_3(x)$ pour tout $x$ de $D_{u_3}$. Pour avoir la contradiction cherch\'ee, il suffit donc de voir que pour tout $x$ de $D_{u_3}$, on a 
$$h'_0({h'_1}^{-1}[h'_2(x)]) = h'_0({h'_1}^{-1}[h'_3(x)]).$$
Posons $H_j := h'_0~\circ~{h'_1}^{-1}~\circ~h'_j$, pour $j=2\mbox{,}~3$, et soit $k\in\omega$.\bigskip

\noindent 1. Pour tout $1\leq i\leq \vert u_0\vert$, $k\not\equiv -1~({J(u_0\lceil i)
\over q_{i-1}})$.\bigskip

 On a alors $H_j(x)(k) = {h'_1}^{-1}[h'_j(x)](k) = h'_1({h'_1}^{-1}[h'_j(x)])(k)$ 
car pour tout $1\leq i\leq \vert u_1\vert $, $k\not\equiv -1~({J(u_1\lceil i)
\over q_{i-1}})$. D'o\`u $H_j(x)(k) = h'_j(x)(k) = x(k)$ car pour tout $1\leq i\leq 
\vert u_j\vert $, $k\not\equiv -1~({J(u_j\lceil i)\over q_{i-1}})$.\bigskip

\noindent 2. Il existe $1\leq i\leq \vert u_0\vert$ maximal tel que $k\equiv -1~({J(u_0\lceil i)
\over q_{i-1}})$, et $q$ tel que $k = {J(u_0\lceil i)\over q_{i-1}}q-1$.\bigskip

\noindent 2.1. L'entier $i$ est maximal sous $\vert u_1\vert $ tel que $k\equiv -1~({J(u_1\lceil i)
\over q_{i-1}})$.
$$\begin{array}{ll}
\mbox{On~a~alors}~ H_j(x)(k)\!\!\!\! & = {h'_1}^{-1}[h'_j(x)](J(u_0\lceil i)q-{J(u_0\lceil i)\over q_{i-1}}-1) \cr 
& = h'_1({h'_1}^{-1}[h'_j(x)])(k) \cr 
& = h'_j(x)(k) \cr 
& = x(J(u_0\lceil i)q-{J(u_0\lceil i)\over q_{i-1}}-1).
\end{array}$$
2.2. L'entier $k$ est de la forme ${J(u_1\lceil \vert u_0\vert+1)\over q_{\vert u_0\vert}}q'-1$ 
(d'o\`u $q = q_{\vert u_0\vert-1}q_{\vert u_0\vert}^{u_1(\vert u_0\vert)}q'$).\bigskip

 On a alors $H_j(x)(k) = {h'_1}^{-1}[h'_j(x)](J(u_0)q-{J(u_0)\over q_{\vert u_0\vert-1}}-1)$. 
Si $1\leq l<\vert u_1\vert $ et $q''$ n'est pas multiple de $q_{l-1}q_l^{u_1(l)}$, on a 
$J(u_0)q-{J(u_0)\over q_{\vert u_0\vert-1}}-1\not= J(u_1\lceil l)q''-{J(u_1\lceil l)\over 
q_{l-1}}-1$ (on simplifie par $q_{\vert u_0\vert-1}^{u_0(\vert u_0\vert-1)}$ si $l>\vert u_0\vert$, 
et par $q_{l-1}^{u_0(l-1)}$ si $l<\vert u_0\vert$ ; si $l=\vert u_0\vert$, on obtient ${q''=q=
q_{\vert u_0\vert-1}q_{\vert u_0\vert}^{u_1(\vert u_0\vert)}q'}$, ce qui est exclus). De m\^eme, on a 
pour tout $q''$ que ${J(u_0)q-{J(u_0)\over q_{\vert u_0\vert-1}}-1\not= J(u_1)q''-
{J(u_1)\over q_{\vert u_1\vert -1}}-1}$. Par ailleurs, 
${J(u_0)q-{J(u_0)\over q_{\vert u_0\vert-1}}-1}={{J(u_1\lceil \vert u_0\vert)\over q_{\vert u_0\vert-1}}(q_{\vert u_0\vert-1}q-1)\! -\! 1}$ ; comme la fonction $h'_1$ est injective, la coordonn\'ee num\'ero ${J(u_0)q-{J(u_0)\over q_{\vert u_0\vert-1}}-1}$ est constante sur $D_{s_1}$, d'o\`u 
$H_2(x) = H_3(x)$.$\hfill\square$

\vfill\eject

\begin{thm} Le th\'eor\`eme 7 devient faux si on suppose seulement $A$ \`a coupes 
verticales d\'enombra-bles.\end{thm} 

\noindent\bf D\'emonstration.\rm\ On raisonne par l'absurde, ce qui fournit un bor\'elien $B_1$. 
Avec ${A = B}$, on voit que $B_1$ a ses coupes horizontales et verticales d\'enombrables. 
Avec ${A = B_1}$, on voit que ${B_1\notin \mbox{pot}(G_\delta )}$. Par le corollaire 8, on obtient 
une situation d'arriv\'ee ${(Z,T,(g_{m,p})_{(m,p)\in\omega^2})}$ et des injections conti-nues 
${i:Z\rightarrow\omega^\omega}$, ${j:T\rightarrow\omega^\omega}$ telles que pour tout $x$ dans 
$G(g)$ on ait ${\overline{g[x]} \cap (i\times j)^{-1}(B_1)_x = g[x]}$. Par le th\'eor\`eme 1, 
on trouve des injections continues ${\tilde u:F\rightarrow G(g)}$ et 
${\tilde v:F\rightarrow T}$ telles que pour tout ${(x,y)\in \bigcup_{n\in\omega} 
\mbox{Gr}(f_n)}$, on ait ${\tilde v(y)\in g[\tilde u(x)]}$, et pour 
$${(x,y)\in (\bigcap_{n\in\omega} D_{f_n}\times F)\! \cap\! 
\overline{\bigcup_{n\in\omega} \mbox{Gr}(f_n)}\setminus 
(\bigcup_{n\in\omega} \mbox{Gr}(f_n))}\mbox{,}$$ 
on ait ${\tilde v(y)\in\overline{g[\tilde u(x)]}\setminus g[\tilde u(x)]}$. Soit 
${(l_s)_{s\in (\omega\setminus\{ 0\} )^{<\omega}}}$ la bonne suite fournie par le lemme 
11 appliqu\'e \`a la suite ${(h_s)_{s\in (\omega\setminus\{ 0\} )^{<\omega}}}$ de 
l'exemple. En appliquant le lemme 10 \`a la bonne suite ${(l_s)_{s\in (\omega\setminus\{ 
0\} )^{<\omega}}}$ et \`a ${G\!  :=\!  \bigcap_{s\in (\omega\setminus\{ 0\} )^{<\omega}} 
D_{l_s}}$, 
on voit que ${A\!  := \! \bigcup_{s\in (\omega\setminus\{ 0\} )^{<\omega}} \mbox{Gr}(l_s\lceil G)
\notin \mbox{pot}(G_\delta )}$. Par suite, on peut trouver des injections continues 
${u:\omega^\omega\rightarrow G}$ et ${v:\omega^\omega\rightarrow 2^\omega}$ telles que 
${\overline{B_1}\cap (u\times v)^{-1}(A) = B_1}$. Posons ${U\! :=\! u\circ i\circ\tilde u}$ 
et ${V\! :=\! v\circ j\circ\tilde v}$ ; $U$ et $V$ sont des hom\'eomorphismes sur leurs 
images (incluses respectivement dans $G$ et $2^\omega$) et si ${(x,y)}\in {\bigcup_{n\in
\omega} \mbox{Gr}(f_n)}$, on a ${(U(x),V(y))\in A}$. De plus, on a 
$$\overline{\bigcup_{n\in\omega} \mbox{Gr}(f_n\lceil G(f))} \cap (U\lceil G(f)\times V)^{-1}(A) 
= \bigcup_{n\in\omega} \mbox{Gr}(f_n\lceil G(f)).$$
$\bullet$ Nous allons montrer un r\'esultat interm\'ediaire. Soient $(s,t)\in 
\bigcup_{n\in\omega} \omega^n\times \omega^{n+1}$, 
$X$ un ouvert non vide de $D_{f_{s,t}}$, et $\tilde s,~\tilde t\in 
(\omega\setminus\{ 0\})^{<\omega}\setminus\{\emptyset\}$ 
telle que pour tout $x\in X$, $U(x)\in V_{\tilde s,\tilde t}$ et 
$(U(x),V(f_{s,t}(x)))\in \mbox{Gr}(l_{\tilde s})$. Alors on peut trouver un ouvert non vide $Y$ 
de $X$, des entiers $m$ et $n$, et $\tilde s',~\tilde t' \in (\omega\setminus\{ 0\})^{<
\omega}$ tels que\bigskip

(a) La suite $\tilde s$ est un d\'ebut strict de $\tilde s'$.\smallskip

(b) L'ouvert $Y$ est inclus dans $D_{f_{s^\frown m,t^\frown n}}$.\smallskip

(c) Pour tout $x$ de $Y$, $U(x)\in V_{\tilde s',\tilde t'}$ et 
$(U(x),V(f_{s^\frown m,t^\frown n}(x)))\in \mbox{Gr}(l_{\tilde s'})$.\smallskip

(d) L'ensemble $l_{\tilde s'} [U[Y]]$ est inclus dans $l_{\tilde s} [U[X]]$.\bigskip

\noindent Soit $O$ (respectivement $P$) un ouvert de $D_{l_{\tilde s}}$ (respectivement 
$2^\omega$) tel que ${U[X] \! =\!  U[F]\! \cap\!  O}$ (respectivement ${V[f_{s,t}[X]] = V[F]\cap P}$). 
Fixons $\alpha\in U[X\cap G(f)]$, $u'\in 2^{<\omega}$ tel que $\alpha\in N_{u'}
\subseteq O$, et, en utilisant le lemme 13,
$$W := \{(x,y)\in N_{v(\tilde s,u',\alpha )}\times P~/~\forall~i<\vert\tilde s\vert~~
y\not= l_{\tilde s\lceil i}(x)~\mbox{et}~d(y,l_{\tilde s}(x)) < \varepsilon (\tilde s,u',\alpha )\}.$$ 
Alors $W$ est ouvert de $2^\omega\times2^\omega$. Comme $U^{-1}(\alpha )\in G(f)$, 
$f_{s^\frown m}(U^{-1}(\alpha ))$ tend vers $f_{s}(U^{-1}(\alpha ))$ quand $m$ tend 
vers l'infini. Comme $U^{-1}(\alpha )\in X\cap U^{-1}(N_{v(\tilde s,u',\alpha )})$, 
on peut trouver un entier $m$ tel que 
${(U^{-1}(\alpha ),f_{s^\frown m}(U^{-1}(\alpha )))\in (U\times V)^{-1}(W)}$. On peut 
trouver $t'\in \omega^{\vert t\vert }$ et $n\in\omega$ tels que 
$f_{s^\frown m}(U^{-1}(\alpha )) = f_{s^\frown m,t'^\frown n}(U^{-1}(\alpha ))$. Comme 
${U^{-1}(\alpha )\in D_{f_{s^\frown m,t'^\frown n}}\subseteq D_{f_{s,t'}}}$, on a $t'=t$ 
puisque ${U^{-1}(\alpha )\in D_{f_{s,t}}}$. Posons  
$$Q := \{x\in X\cap U^{-1}(N_{v(\tilde s,u',\alpha )})\cap D_{f_{s^\frown m,t^
\frown n}}~/~(U(x),V(f_{s^\frown m,t^\frown n}(x)))\in W\}.$$

\vfill\eject

 Alors $Q\subseteq \bigcup_{w\in (\omega\setminus\{ 0\})^{<\omega}} \{x\in Q~/~
(U(x),V(f_{s^\frown m,t^\frown n}(x)))\in \mbox{Gr}(l_w)\}$. Par le th\'eor\`eme de Baire, on 
peut donc trouver un ouvert non vide $Y$ de $Q$ et $\tilde s'\in 
(\omega\setminus\{ 0\})^{<
\omega}$ tels que pour tout $x$ de $Y$, $(U(x),V(f_{s^\frown m,t^\frown n}(x)))\in 
\mbox{Gr}(l_{\tilde s'})$. Par le th\'eor\`eme de Baire encore, on peut supposer qu'il existe 
$\tilde t'\in (\omega\setminus\{ 0\})^{\vert\tilde s'\vert+1}$ telle que pour tout $x$ de 
$Y$, $U(x)\in V_{\tilde s',\tilde t'}$.\bigskip

 Les conditions (b) et (c) sont clairement v\'erifi\'ees. Si $x$ est dans $Y$, on a l'\'egalit\'e 
$${l_{\tilde s'}(U(x)) = V(f_{s^\frown m,t^\frown n}(x))}.$$
Comme $x$ est dans $Q$, $(U(x),V(f_{s^\frown m,t^\frown n}(x)))$ est dans $W$, donc 
${V(f_{s^\frown m,t^\frown n}(x))\!  \in \! V[f_{s,t}[X]]}$. Donc on trouve $y\in X$ tel 
que ${V(f_{s^\frown m,t^\frown n}(x))\!  =\!  
l_{\tilde s}(U(y))}$. D'o\`u la condition (d). De plus, ${V(f_{s^\frown m,t^\frown n}(x))
\not= l_{\tilde s\lceil i}(U(x))}$ si $i<\vert\tilde s\vert$, et on a \'egalement 
$d(l_{\tilde s}(U(x)),l_{\tilde s'}(U(x))) < \varepsilon (\tilde s,u',\alpha )$. Donc $\tilde s$ 
est un d\'ebut de $\tilde s'$ puisque ${U(x)\in N_{v(\tilde s,u',\alpha )}\cap G}$. Si 
$\tilde s = \tilde s'$, on a successivement 
$${V(f_{s^\frown m,t^\frown n}(x)) = l_{\tilde s'}(U(x)) = 
l_{\tilde s}(U(x)) = V(f_{s,t}(x))}\mbox{,}$$ 
d'o\`u l'\'egalit\'e ${f_{s^\frown m,t^\frown n}(x) = f_{s,t}(x)}$, qui est absurde.\bigskip

\noindent $\bullet$ Revenons \`a la preuve du th\'eor\`eme. On peut trouver $(c,d)\in\bigcup_{n\in\omega} 
\omega^{n}\times \omega^{n+1}$ et un ouvert non vide $R$ de $D_{f_{c,d}}$ tels que 
pour tout $x$ de $R$ on ait $(U(x),V(f_{c,d}(x)))\notin \mbox{Gr}(l_\emptyset)$. Sinon 
$$H := \bigcap_{n\in\omega} \{x\in G(f)~/~(U(x),V(f_n(x)))\in \mbox{Gr}(l_\emptyset)\}$$ 
serait $G_\delta$ dense de $F$ et on aurait $\bigcup_{n\in\omega} \mbox{Gr}(f_n\lceil H) = 
\overline{\bigcup_{n\in\omega} \mbox{Gr}(f_n\lceil H)}\cap (U\lceil H\times V)^{-1}
(\mbox{Gr}(l_\emptyset\lceil G))$, donc $\bigcup_{n\in\omega} \mbox{Gr}(f_n\lceil H)$ serait 
$\mbox{pot}(\bormone)$ non $\mbox{pot}(\bormtwo)$, ce qui est absurde.\bigskip

 On a l'inclusion 
$$R \subseteq \bigcup_{w\in (\omega\setminus\{ 0\})^{<\omega}} 
\{x\in R~/~(U(x),V(f_{c,d}(x)))\in \mbox{Gr}(l_w)\}.$$
Par le th\'eor\`eme de Baire comme ci-dessus, on trouve un ouvert non vide $X$ de 
$R$ et $\tilde s,\tilde t$ dans 
$(\omega\setminus\{ 0\})^{<\omega}\setminus\{\emptyset\}$ 
tels que pour tout $x$ de $X$, on ait $U(x)\in V_{\tilde s,\tilde t}$ et 
${(U(x),V(f_{c,d}(x)))\in \mbox{Gr}(l_{\tilde s})}$. On applique le point pr\'ec\'edent \`a 
$c$, $d$, et $X$, ce qui fournit $Y_0$, $m_0$, $n_0$, et $s_0,~t_0$. On 
applique ensuite le point pr\'ec\'edent \`a $c^\frown m_0$, $d^\frown n_0$ et $Y_0$, ce qui fournit 
$Y_1$, $m_1$, $n_1$, et $s_1,~t_1$. On applique ensuite le point pr\'ec\'edent \`a 
$c^\frown m_0^\frown m_1$, $d^\frown n_0^\frown n_1$ et $Y_1$, ce qui fournit $Y_2$, $m_2$, 
$n_2$, et $s_2,~t_2$. On applique enfin le point pr\'ec\'edent \`a 
$c^\frown m_0^\frown m_1^\frown m_2$, $d^\frown n_0^\frown n_1^\frown n_2$ et $Y_2$, ce qui 
fournit $Y_3$, $m_3$, $n_3$, et $s_3,~t_3$. On a $\tilde s \prec_{\not=} s_0$ et $s_i  
\prec_{\not=} s_{i+1}$ si $i\in 2$, donc les $s_i$ sont non vides. 
Posons $D_{s_i} := U[Y_i]$. Avec les notations du lemme 14, on a $\mbox{Im}(l'_{i+1}) = 
l_{s_{i+1}}[U[Y_{i+1}]] \subseteq l_{s_{i}}[U[Y_{i}]] = \mbox{Im}(l'_i)$. Pour tout $x$ de $Y_i$, on a 
$V(f_{m_0^\frown ...^\frown m_i,0^\frown n_0^\frown ...^\frown n_i}(x))\! =\! l_{s_i}(U(x))$. Par suite, si 
$y\! =\! U(x)$ et $y'\! =\! U(x')$ avec $x,~x'
\in Y_i$, $l'_i(y) = l'_i(y')$ entra\^\i ne successivement que 
$V(f_{m_0^\frown ...^\frown m_i,0^\frown n_0^\frown ...^
\frown n_i}(x))$ vaut $V(f_{m_0^\frown ...^\frown m_i,0^\frown n_0^\frown ...^
\frown n_i}(x'))$, que $f_{m_0^\frown ...^\frown m_i,0^\frown n_0^\frown ...^
\frown n_i}(x) = f_{m_0^\frown ...^\frown m_i,0^\frown n_0^\frown ...^
\frown n_i}(x')$, puis que $x = x'$ et $y = y'$. D'o\`u l'injectivit\'e de $l'_i$. Le lemme 14 
peut donc s'appliquer et donne la contradiction cherch\'ee.$\hfill\square$

\vfill\eject

\noindent\bf {\Large 5 R\'ef\'erences.}\rm\bigskip

\noindent [HKL]\ \ L. A. Harrington, A. S. Kechris et A. Louveau,~\it A Glimm-Effros 
dichotomy for Borel equivalence relations,~\rm J. Amer. Math. Soc.~3 (1990), 903-928

\noindent [Ke]\ \ A. S. Kechris,~\it Classical Descriptive Set Theory,~\rm 
Springer-Verlag, 1995

\noindent [Ku]\ \  K. Kuratowski,~\it Topology,~\rm Vol. 1, Academic Press, 1966

\noindent [Le1]\ \ D. Lecomte,~\it Classes de Wadge potentielles et th\'eor\`emes d'uniformisation 
partielle,~\rm Fund. Math.~143 (1993), 231-258

\noindent [Le2]\ \ D. Lecomte,~\it Uniformisations partielles et crit\`eres \`a la Hurewicz dans le plan,~\rm Trans. A.M.S. ~347, 11 (1995), 4433-4460

\noindent [Le3]\ \ D. Lecomte,~\it Tests \`a la Hurewicz dans le plan,~\rm Fund. Math.~156 (1998), 
131-165

\noindent [Lo1]\ \ A. Louveau,~\it Ensembles analytiques et bor\'eliens dans les espaces produit,~\rm 
Ast\'erisque (S. M. F.) 78 (1980)

\noindent [Lo2]\ \ A. Louveau,~\it A separation theorem for $\Ana$ sets,~\rm Trans. A. M. S.~260 (1980), 363-378

\noindent [Lo3]\ \ A. Louveau,~\it Livre \`a para\^\i tre\rm

\noindent [Lo-SR]\ \ A. Louveau and J. Saint Raymond,~\it Borel classes and closed games : Wadge-type and Hurewicz-type results,~\rm Trans. A. M. S.~304 (1987), 431-467

\noindent [Mo]\ \ Y. N. Moschovakis,~\it Descriptive set theory,~\rm North-Holland, 1980

\noindent [SR]\ \ J. Saint Raymond,~\it La structure bor\'elienne d'Effros est-elle standard ?,~\rm 
Fund. Math.~100 (1978), 201-210

\end{document}